\documentclass[11pt]{amsart}
\usepackage[top=2.2cm,left=2.8cm,right=2.8cm,bottom=2.8cm]{geometry}
\usepackage[hyperfootnotes=false]{hyperref}
\usepackage{amsmath,amssymb,mathtools}
\usepackage{graphicx,booktabs,multirow}
\usepackage[shortlabels]{enumitem}
\usepackage{xcolor}
\usepackage{algorithm}
\usepackage{algpseudocode}
\hypersetup{colorlinks=true,linkcolor=blue!50!black,citecolor=blue!50!black,urlcolor=blue!50!black}
\numberwithin{equation}{section}

\algrenewcommand\algorithmicrequire{\textbf{Input:}}
\algrenewcommand\algorithmicensure{\textbf{Output:}}

\newcommand{\T}{\mathbb{T}^d}
\newcommand{\dd}{\,\mathrm{d}}
\newcommand{\eps}{\varepsilon}
\newcommand{\F}{\mathcal{F}}
\newcommand{\Tb}{\mathcal{T}_{\rm b}}
\newcommand{\Tf}{\mathcal{T}_{\rm f}}

\newcommand{\I}{I}

\usepackage{microtype}
\usepackage{placeins}
\setlist[enumerate]{itemsep=2pt,parsep=0pt,topsep=4pt}
\usepackage[labelfont=bf,labelsep=period,font=small]{caption}

\theoremstyle{plain}
\newtheorem{theorem}{Theorem}[section]
\newtheorem{lemma}[theorem]{Lemma}
\newtheorem{proposition}[theorem]{Proposition}
\newtheorem{corollary}[theorem]{Corollary}
\theoremstyle{definition}
\newtheorem{definition}[theorem]{Definition}
\newtheorem{remark}[theorem]{Remark}
\numberwithin{equation}{section}

\title[Newton iterations and SL schemes for MFGC]{Newton iterations and semi-Lagrangian schemes for mean field games of controls}
\author{Ahmad Zorkot}
\thanks{Faculty of Mathematics, University of Vienna, Oskar-Morgenstern-Platz 1, 1090 Vienna, Austria. \emph{E-mail address}: \href{mailto:ahmad.zorkot@univie.ac.at}{ahmad.zorkot@univie.ac.at}}

\begin{document}
\begin{abstract}
We study Newton methods for second-order mean field games of controls (MFGC), in which agents interact through local averages of their controls. We prove local quadratic convergence of full Newton iterations, and of Newton iterations on reduced unknowns, near stable classical solutions. Stability holds under weak control interaction and under an infinitesimal monotonicity condition of Lasry--Lions type. We propose a semi-Lagrangian (SL) scheme for second-order MFGC, in which the average-control relation is solved by a convergent Chebyshev iteration, and two matrix-free solvers: a full Newton--SL method with block Gauss--Seidel preconditioning, whose spectrum we characterize in the absence of control coupling, and a reduced Jacobian-free Newton--GMRES method whose density iterates are nonnegative and mass-preserving by construction. We also construct an exact, time-dependent, coupled MFGC solution for use as a benchmark. In the experiments, the reduced solver converges quadratically in its final iterations, also in a strong-alignment test in which full Newton updates lose positivity. On the benchmark, the SL scheme and a finite-difference scheme both converge at first order with errors of comparable size, and our SL implementation is faster.


\end{abstract}

\maketitle

\section{Introduction}

Mean field games (MFG) provide a framework for describing Nash equilibria in differential games with a large number of indistinguishable agents \cite{LL07,HMC06}. In many applications, the cost of an agent depends on the controls chosen by the others: the price of a good depends on the aggregate production \cite{CL18,GV16}, a trader is affected by the average trading rate of the market \cite{CL18}, and in crowd motion a pedestrian may wish to move in the same direction as her neighbors \cite{AK21}. These models are called mean field games of controls (MFGC), or extended mean field games \cite{GPV14,GV16}, see also \cite{CarLac15}. An MFGC system is described by a backward Hamilton--Jacobi--Bellman (HJB) equation, a forward Fokker--Planck (FP) equation and a fixed-point relation for the law of the controls. Existence and uniqueness of classical solutions were obtained by Kobeissi \cite{Kob22} under a structural condition which makes the fixed-point relation a contraction, and under a monotonicity condition of Lasry--Lions type \cite{Kob20}. For a class of potential MFGC, existence and regularity were studied in \cite{BHP21}.

From the numerical point of view, MFGC were studied by Achdou and Kobeissi \cite{AK21}, who extended to MFGC the finite difference schemes of \cite{AC10}, with a discrete FP equation defined as the adjoint of the linearized discrete HJB equation. Their discrete system is solved by a Newton method on reduced unknowns, the coupling term and the average drift, while the value function and the density are obtained by implicit time marching; in their crowd motion example with strong alignment, they use continuation in the viscosity parameter, and they observe several distinct solutions. Meynard \cite{Mey26} proposed extragradient methods for MFGC, which converge exponentially under strong monotonicity assumptions. Newton methods for MFG have been studied at the discrete level \cite{AC10,AP12} and at the continuous level: Camilli and Tang \cite{CT} proved the local quadratic convergence of the Newton iterations for second-order MFG in the monotone case, and Berry, Ley and Silva \cite{BLS} established the corresponding result for stationary MFG around stable solutions, a notion introduced by Briani and Cardaliaguet \cite{BC18}. Carlini and Zorkot \cite{CZ} discretized the Newton iterations for MFG by semi-Lagrangian (SL) schemes, and solved each linear system by a block Gauss--Seidel method, in which the matrices are never assembled. SL schemes for MFG were introduced by Carlini and Silva \cite{CS14,CS15}; they are explicit, do not require a parabolic CFL condition, and their discrete FP equation, defined with the transpose of the SL transport operator, preserves positivity and mass. A Lagrange--Galerkin scheme for first-order MFG, also based on characteristics, was analyzed in \cite{CSZ24}. A fully discrete SL scheme for a first-order price formation MFG, in which the price is determined by the aggregate control, was proposed in \cite{AG25}. To the best of our knowledge, the convergence of Newton iterations for MFGC has not been analyzed at the continuous level, and no SL scheme has been proposed for second-order MFGC.

The aim of the present work is to address this gap, and its main focus is the solution of MFGC by Newton iterations, both at the continuous level and for an SL discretization. We first prove the local quadratic convergence of the Newton iterations for MFGC around any stable solution, in any space dimension (Theorem~\ref{thm:main}). The proof adapts to the time-dependent MFGC system the strategy of \cite{BLS} for stationary MFG; the new difficulty is the relation defining the average control, which is not compact, and which we eliminate by a Schur complement, using the contraction condition of \cite{Kob22}. A classical solution of the problem without control interaction persists as a stable solution under sufficiently weak control interactions (Corollary~\ref{cor:small}), and a solution is stable under an infinitesimal form of the monotonicity condition of \cite{Kob20} (Proposition~\ref{prop:monotone}). We also prove the local quadratic convergence of Newton iterations on reduced unknowns (Theorem~\ref{thm:reduced}), the continuous counterpart of the strategy of \cite{AK21}, along which the density is nonnegative and of unit mass.

We then propose an SL scheme for MFGC, in which the relation defining the average control is solved by a Chebyshev iteration, convergent when the density is nonnegative (Lemma~\ref{lem:Zdiscrete}), and two matrix-free Newton solvers for it. The first one is the Newton--SL method of \cite{CZ}, extended to MFGC; we show that without interaction through the controls its block Gauss--Seidel sweep has a real nonpositive spectrum (Proposition~\ref{prop:bgs}), while for alignment the dominant eigenvalue can be larger than one (Section~\ref{sec:test1}). The second one is an SL realization of the reduced iterations, solved by a Jacobian-free Newton--GMRES method, while preserving nonnegativity and the initial discrete mass of the reconstructed density at every iteration. We compare the two reduced discretizations on an independently verifiable time-dependent coupled MFGC benchmark. The reduced unknowns are taken from \cite{AK21}, and Newton--Krylov methods and nonlinear elimination are classical \cite{KK04,LRW96}; we do not claim them as new.

To the best of our knowledge, the following are new: (i) the convergence analysis of Newton iterations for MFGC at the continuous level; (ii) an SL scheme for second-order MFGC; (iii) a matrix-free reduced Newton solver for it, along which the density is nonnegative and of unit mass. We also construct an exact, time-dependent, coupled MFGC solution, in which only the HJB equation contains a prescribed term, and use it to measure discretization errors.

The numerical experiments distinguish solver robustness from discretization accuracy. Full Newton--SL converges on a strong-alignment test with suitable initialization, whereas direct Newton updates can violate density positivity; the reduced formulation is more robust in that case. Reduced Newton--SL, which differentiates the discrete residual, converges quadratically in its final iterations, whereas full Newton--SL, which discretizes the continuous linearization and is therefore not the exact Newton method for the discrete system, converges linearly. Reduced Newton--SL also converges in a few iterations under spatial refinement on a smooth coupled problem. Independent errors in the value function, density and averaged control are measured using an exact, time-dependent coupled MFGC with an unforced FP equation and prescribed running potential. On this benchmark, reduced Newton--SL and the reduced finite-difference Newton method both converge at first order under refinement with $\Delta t\propto h$, with errors of comparable size, and our SL implementation is faster. In the crowd experiment, computed density patterns agree qualitatively with those reported in \cite{AK21}; this is not an independent error estimate.

The convergence results for Newton iterations are local and proved on the torus, whereas the crowd example is on a bounded square with reflecting boundaries. We do not prove convergence of the fully discrete, non-monotone SL scheme, and do not infer continuum multiplicity from the computed crowd branches. The exact time-dependent test uses a prescribed running potential $V(t,x)$ and is a validation of the discretization, not a further convergence theorem.

The paper is organized as follows. Section~\ref{sec:model} presents the model and the assumptions. Section~\ref{sec:theory} contains the convergence analysis of the Newton iterations at the continuous level. Section~\ref{sec:sl} introduces the SL scheme, and Section~\ref{sec:solvers} the two Newton solvers. Section~\ref{sec:numerics} reports the numerical experiments, and Section~\ref{sec:conclusions} contains the conclusions. The exact benchmark solution is constructed in Appendix~\ref{app:exactMFGC}.

\section{The model}\label{sec:model}

Let $Q=(0,T)\times\T$, where $\T$ is the $d$-dimensional torus. In this work, we consider the second-order MFGC system
\begin{equation}\label{eq:mfgc}
\begin{cases}
-\partial_t u-\nu\Delta u+H(x,Du,Z)=f(x,m) & \text{in } Q,\\[2pt]
\partial_t m-\nu\Delta m-\operatorname{div}\big(m\,D_pH(x,Du,Z)\big)=0 & \text{in } Q,\\[2pt]
Z=\Phi(u,m,Z) & \text{in } Q,\\[2pt]
u(T,x)=G(x),\quad m(0,x)=m_0(x) & \text{in } \T,
\end{cases}
\end{equation}
where $\nu>0$ is the diffusion coefficient, $u$ is the value function of a representative agent, $m$ is the density of the agents, and $Z$ is the local average of the optimal controls,
\begin{equation}\label{eq:Phi}
\Phi(u,m,Z)(t,x)=\frac{\int_{\T}K(x,y)\big(-D_pH(y,Du(t,y),Z(t,y))\big)\,m(t,y)\dd y}{\int_{\T}K(x,y)\,m(t,y)\dd y+\eps},\qquad \eps>0 .
\end{equation}
Here $K\ge0$ is an interaction kernel. A representative agent controls its state through $\dd X_t=\alpha_t\dd t+\sqrt{2\nu}\dd W_t$ and minimizes
\[
\mathbb E\Big[\int_t^TL(X_s,\alpha_s,Z(s,X_s))+f(X_s,m(s,X_s))\dd s+G(X_T)\Big],
\]
where the Lagrangian $L$ and the Hamiltonian $H$ are related by
\begin{equation}\label{eq:legendreH}
H(x,p,Z)=\sup_{\alpha\in\mathbb R^d}\big\{-\alpha\cdot p-L(x,\alpha,Z)\big\}.
\end{equation}

Our main example is the crowd motion model of \cite{AK21}, in which an agent pays for deviating from the average control of her neighbors:
\begin{equation}\label{eq:align}
L(x,\alpha,Z)=\tfrac12|\alpha-\lambda Z|^2,\qquad H(x,p,Z)=\tfrac12|p|^2-\lambda Z\cdot p,\qquad \alpha^*=\lambda Z-Du,
\end{equation}
where $|\lambda|<1$. Positive $\lambda$ describes alignment, negative $\lambda$ anti-alignment, and $\lambda=0$ reduces to a classical MFG.

\subsection{Assumptions}
We make the following assumptions:
\begin{enumerate}[(i)]
\item $\nu>0$; $m_0\in C^{2+\beta}(\T)$, $m_0\ge0$, $\int_{\T}m_0=1$; $G\in C^{2+\beta}(\T)$ for some $\beta\in(0,1)$.
\item $H\in C^3$, $D^2_{pp}H>0$, and the third derivatives of $H$ are locally Lipschitz.
\item $f\in C^2(\T\times\mathbb R)$ and $\partial^2_mf$ is locally Lipschitz.
\item $K\in C^\infty(\T\times\T)$, $K\ge0$.
\item There is $\lambda_0\in[0,1)$ such that $|D_pH(x,p,Z_1)-D_pH(x,p,Z_2)|\le\lambda_0|Z_1-Z_2|$ for all $x,p,Z_1,Z_2$.
\end{enumerate}

Note that assumption (v) is the structural condition of Kobeissi \cite{Kob22}, which makes the relation defining $Z$ a contraction; for the alignment model \eqref{eq:align}, $\lambda_0=|\lambda|$.

\section{Newton iterations at the continuous level}\label{sec:theory}

\subsection{The Newton iterations}
Assuming sufficient smoothness of the data, we write \eqref{eq:mfgc} as $\F(u,m,Z)=0$, where $\F$ collects the residuals of the HJB and FP equations, with their terminal and initial conditions, and of \eqref{eq:Phi}; it is defined precisely, in mild form, in the proof of Theorem~\ref{thm:main}. Given $x^{n-1}=(u^{n-1},m^{n-1},Z^{n-1})$, Newton's method defines $x^n=(u^n,m^n,Z^n)$ by
\[
D\F(x^{n-1})(x^n-x^{n-1})=-\F(x^{n-1}).
\]
We denote by $H^{n-1}$, $H^{n-1}_Z$, $H^{n-1}_{pp}$, $H^{n-1}_{pZ}$ the values of $H$ and of its derivatives at $(x,Du^{n-1},Z^{n-1})$, and by $\partial\Phi^{n-1}$ the partial derivatives of $\Phi$ at $(u^{n-1},m^{n-1},Z^{n-1})$. Setting $q^n=D_pH(x,Du^{n-1},Z^{n-1})$, the Newton step amounts to substituting in \eqref{eq:mfgc} the first-order expansions of $H$, $f$ and $\Phi$ around $x^{n-1}$, which gives the linearized system
{\small
\begin{equation}\label{eq:newton}
\begin{cases}
-\partial_t u^n-\nu\Delta u^n+q^n\cdot Du^n+H_Z^{n-1}\cdot Z^n-f'(m^{n-1})\,m^n=R^{n-1},\\[2pt]
\partial_t m^n-\nu\Delta m^n-\operatorname{div}(m^nq^n)-\operatorname{div}\big(m^{n-1}[H_{pp}^{n-1}Du^n+H_{pZ}^{n-1}Z^n]\big)=S^{n-1},\\[2pt]
Z^n-\partial_u\Phi^{n-1}u^n-\partial_m\Phi^{n-1}m^n-\partial_Z\Phi^{n-1}Z^n=P^{n-1},\\[2pt]
u^n(T)=G,\quad m^n(0)=m_0,
\end{cases}
\end{equation}}
with
\[
\begin{aligned}
R^{n-1}&=f(m^{n-1})-f'(m^{n-1})m^{n-1}-H^{n-1}+q^n\cdot Du^{n-1}+H^{n-1}_Z\cdot Z^{n-1},\\
S^{n-1}&=-\operatorname{div}\big(m^{n-1}[H_{pp}^{n-1}Du^{n-1}+H_{pZ}^{n-1}Z^{n-1}]\big),\\
P^{n-1}&=\Phi^{n-1}-\partial_u\Phi^{n-1}u^{n-1}-\partial_m\Phi^{n-1}m^{n-1}-\partial_Z\Phi^{n-1}Z^{n-1}.
\end{aligned}
\]
For $H=H(x,p)$, the first two equations are the Newton iterations for MFG studied in \cite{CT,CZ}. Each Newton step is a coupled forward--backward linear system, which includes a linear relation for $Z^n$. For the alignment model \eqref{eq:align}, $q^n=Du^{n-1}-\lambda Z^{n-1}$, $H^{n-1}_Z=-\lambda Du^{n-1}$, $H_{pp}=I$ and $H_{pZ}=-\lambda I$. We write $f'(m)=\partial_mf(x,m)$, suppressing the spatial argument.

\subsection{Quadratic convergence}
Let
\[
\begin{gathered}X=C^{1,0}(\bar Q)\times C^0(\bar Q)\times C^0(\bar Q;\mathbb R^d),\\ \|(u,m,Z)\|_X=\|u\|_{\infty}+\|Du\|_\infty+\|m\|_\infty+\|Z\|_\infty ,\end{gathered}
\]
where $C^{1,0}(\bar Q)$ is the space of functions which are continuous on $\bar Q$ together with their spatial gradient; this is the norm in which the Newton errors are measured in \cite{CT}.

Let $\Gamma_\tau$, $\tau>0$, be the heat kernel of $\partial_t-\nu\Delta$ on $\T$, and define, for $g\in C^0(\bar Q)$ and $b\in C^0(\bar Q;\mathbb R^d)$,
\[
\Tb[g](t)=\int_t^T\Gamma_{s-t}*g(s)\dd s,\qquad \Tf[\operatorname{div}b](t)=\int_0^t\nabla\Gamma_{t-s}*b(s)\dd s,
\]
the mild solutions of $-\partial_tv-\nu\Delta v=g$, $v(T)=0$, and $\partial_t\rho-\nu\Delta\rho=\operatorname{div}b$, $\rho(0)=0$; here $*$ is the convolution on $\T$, and $\nabla\Gamma_{\tau}*b=\sum_i\partial_i\Gamma_\tau*b_i$, which equals $\Gamma_\tau*\operatorname{div}b$ for smooth $b$. Let $\bar u$, $\bar m$ solve the homogeneous equations with data $G$, $m_0$.

\begin{lemma}\label{lem:heat}
There is $C>0$, depending on $\nu$, $T$ and $d$, such that for $g\in C^0(\bar Q)$, $b\in C^0(\bar Q;\mathbb R^d)$ and $t\in[0,T]$:
\begin{enumerate}[(a)]
\item the following Volterra-type bounds hold:
\[
\begin{gathered}\|D\,\Tb[g](t)\|_\infty\le C\int_t^T\frac{\|g(s)\|_\infty}{(s-t)^{1/2}}\dd s,\\ \|\Tf[\operatorname{div}b](t)\|_\infty\le C\int_0^t\frac{\|b(s)\|_\infty}{(t-s)^{1/2}}\dd s;\end{gathered}
\]
\item $\Tb[g]$, $D\,\Tb[g]$ and $\Tf[\operatorname{div}b]$ are H\"older continuous on $\bar Q$, with exponent $\frac12$ in $t$ and any $\beta<1$ in $x$, with H\"older seminorms bounded by $C_\beta\|g\|_\infty$ and $C_\beta\|b\|_\infty$;
\item $\int_{\T}\Tf[\operatorname{div}b](t)\dd x=0$ for all $t$.
\end{enumerate}
Consequently, $\Tb:C^0(\bar Q)\to C^{1,0}(\bar Q)$ and $\Tf\circ\operatorname{div}:C^0(\bar Q;\mathbb R^d)\to C^0(\bar Q)$ are compact.
\end{lemma}
\begin{proof}
We use the classical bounds \[\|D^k\Gamma_\tau\|_{L^1(\T)}\le C\tau^{-k/2},\ k=1,2,\qquad \|\partial_\tau D\Gamma_\tau\|_{L^1(\T)}\le C\tau^{-3/2},\qquad 0<\tau\le T.\] Part (a) follows from Young's inequality. For (b), consider $\rho=\Tf[\operatorname{div}b]$; the argument for $D\,\Tb[g]$ is the same, and $\Tb[g]$ is Lipschitz in $x$ and H\"older in $t$ by the same bounds. For $y\in\mathbb R^d$, $\|\nabla\Gamma_\tau(\cdot+y)-\nabla\Gamma_\tau\|_{L^1}\le C\min(\tau^{-1/2},|y|\tau^{-1})$, so that, for $|y|\le\sqrt d$ and with $a=\min\{|y|^2,T\}$,
\[
\begin{aligned}
\|\rho(t,\cdot+y)-\rho(t)\|_\infty&\le C\|b\|_\infty\Big(\int_0^{a}\tau^{-1/2}\dd\tau+\int_{a}^{T}|y|\tau^{-1}\dd\tau\Big)\\
&\le C\|b\|_\infty|y|\big(1+\log^+(T/|y|^2)\big)\le C_\beta\|b\|_\infty|y|^\beta .
\end{aligned}
\]
For $0\le t<t'\le T$ and $\delta=t'-t$, $\rho(t')-\rho(t)=\int_t^{t'}\nabla\Gamma_{t'-s}*b(s)\dd s+\int_0^t(\nabla\Gamma_{t-s+\delta}-\nabla\Gamma_{t-s})*b(s)\dd s$. The first term is bounded by $C\|b\|_\infty\delta^{1/2}$ and, since $\|\nabla\Gamma_{\tau+\delta}-\nabla\Gamma_\tau\|_{L^1}\le C\min(\tau^{-1/2},\delta\tau^{-3/2})$, so is the second. Part (c) follows from $\int_{\T}\nabla\Gamma_\tau*b\dd x=0$. Compactness follows from (b) and the Arzel\`a--Ascoli theorem: for $g$ and $b$ in bounded sets of $C^0$, the functions $\Tb[g]$ and $D\,\Tb[g]$ (hence the compactness into $C^{1,0}(\bar Q)$) and $\Tf[\operatorname{div}b]$ are uniformly bounded and equicontinuous on $\bar Q$.
\end{proof}

We will also use the Gronwall lemma for weakly singular kernels (see \cite[Lemma 7.1.1]{Hen81}): if $\phi\in C^0([0,T])$ is nonnegative, $a\ge0$, and $\phi(t)\le a+C\int_0^t(t-s)^{-1/2}\phi(s)\dd s$ for all $t$, or $\phi(t)\le a+C\int_t^T(s-t)^{-1/2}\phi(s)\dd s$ for all $t$, then $\phi\le C'a$, with $C'$ depending only on $C$ and $T$; in particular $\phi=0$ if $a=0$. Together with Lemma~\ref{lem:heat}(a), it shows that, for $b\in C^0(\bar Q;\mathbb R^d)$, the equations $\rho=\Tf[\operatorname{div}(\rho b)]$ in $C^0(\bar Q)$ and $v=-\Tb[b\cdot Dv]$ in $C^{1,0}(\bar Q)$ have only the trivial solution (apply it to $\phi(t)=\|\rho(t)\|_\infty$ and $\phi(t)=\|Dv(t)\|_\infty$; then $v=0$ since $v=-\Tb[b\cdot Dv]$).

\begin{definition}\label{def:stable}
A classical solution $x^*=(u^*,m^*,Z^*)$ of \eqref{eq:mfgc} is called \emph{stable} if $(0,0,0)$ is the only solution in $X$ of the linearized system
\begin{equation}\label{eq:lin}
\begin{cases}
-\partial_tv-\nu\Delta v+H_p^*\cdot Dv+H_Z^*\cdot\zeta=f'(m^*)\rho,\\[2pt]
\partial_t\rho-\nu\Delta\rho-\operatorname{div}(\rho H_p^*)-\operatorname{div}\big(m^*[H_{pp}^*Dv+H_{pZ}^*\zeta]\big)=0,\\[2pt]
\zeta=\partial_u\Phi^*v+\partial_m\Phi^*\rho+\partial_Z\Phi^*\zeta,\qquad v(T)=0,\quad \rho(0)=0,
\end{cases}
\end{equation}
where $*$ denotes evaluation at $x^*$, understood in the mild sense: $v=\Tb[f'(m^*)\rho-H^*_p\cdot Dv-H^*_Z\cdot\zeta]$ and $\rho=\Tf[\operatorname{div}(\rho H^*_p+m^*(H^*_{pp}Dv+H^*_{pZ}\zeta))]$.
\end{definition}
This is the extension to MFGC of the stability notion used for MFG in \cite{BC18,BLS}. Stability does not require monotonicity.

Mild solutions $(v,\rho,\zeta)\in X$ of \eqref{eq:lin} are weak solutions with the regularity needed for the energy identity below. Set
\begin{equation}\label{eq:Fg}
F=\rho H_p^*+m^*\big(H_{pp}^*Dv+H_{pZ}^*\zeta\big),\qquad g=f'(m^*)\rho-H_p^*\cdot Dv-H_Z^*\cdot\zeta,
\end{equation}
which are continuous on $\bar Q$, so that $-\partial_tv-\nu\Delta v=g$ and $\partial_t\rho-\nu\Delta\rho=\operatorname{div}F$. Since $v\in C^{1,0}(\bar Q)$, we have $v\in L^2(0,T;H^1(\T))$. Moreover, $\rho=\Tf[\operatorname{div}F]$ is the weak solution of its equation with $\rho(0)=0$ (approximate $F$ uniformly by smooth fields and use Lemma~\ref{lem:heat}(a)), and the energy estimate gives $\rho\in L^2(0,T;H^1(\T))$ \cite[Chapter III]{LSU}. The two equations then give $\partial_tv,\partial_t\rho\in L^2(0,T;H^{-1}(\T))$. Consequently, $t\mapsto\int_{\T}\rho(t)v(t)\dd x$ is absolutely continuous on $[0,T]$, with derivative $\langle\partial_t\rho,v\rangle+\langle\partial_tv,\rho\rangle$.

\begin{theorem}\label{thm:main}
Assume (i)--(v) and let $x^*=(u^*,m^*,Z^*)$ be a stable classical solution of \eqref{eq:mfgc}. Then there exist $\eta>0$ and $C>0$ such that, if $\|x^0-x^*\|_X\le\eta$, the Newton iterations \eqref{eq:newton} are well defined for all $n\ge1$, i.e.\ each linear system \eqref{eq:newton} has a unique solution $x^n=(u^n,m^n,Z^n)\in X$, and
\[
\|x^n-x^*\|_X\le C\,\|x^{n-1}-x^*\|_X^2\qquad\text{for all } n\ge1 .
\]
In particular, $x^n\to x^*$ quadratically in $X$.
\end{theorem}

\begin{proof}
\emph{Step 1 (operator formulation).} In mild form, the map $\F$ of Section~\ref{sec:theory} is $\F:\mathcal O\subset X\to X$,
\[
\F(u,m,Z)=\begin{pmatrix}u-\bar u-\Tb\big[f(\cdot,m)-H(\cdot,Du,Z)\big]\\ m-\bar m-\Tf\big[\operatorname{div}\big(mD_pH(\cdot,Du,Z)\big)\big]\\ Z-\Phi(u,m,Z)\end{pmatrix},
\]
on the open set $\mathcal O=\{(u,m,Z)\in X:\ \min_{(t,x)\in\bar Q}\{\int K(x,y)m(t,y)\dd y+\eps\}>0\}$, which contains $x^*$ since $m^*\ge0$. For $(u,m,Z)\in X$, the functions $f(\cdot,m)$, $H(\cdot,Du,Z)$ and $mD_pH(\cdot,Du,Z)$ are continuous on $\bar Q$, so that $\F$ maps $\mathcal O$ into $X$. Under assumptions (ii)--(iii), the superposition operators $m\mapsto f(\cdot,m)$ and $(p,Z)\mapsto H(\cdot,p,Z)$, $D_pH(\cdot,p,Z)$ are of class $C^2$ on spaces of continuous functions \cite[Chapter 1]{AP93}, hence of class $C^1$ with locally Lipschitz derivatives. By Lemma~\ref{lem:heat}, $\Tb$ and $\Tf\circ\operatorname{div}$ are bounded from $C^0$ into $C^{1,0}(\bar Q)$ and $C^0(\bar Q)$. Since $K$ is smooth, $\phi\mapsto\int K(\cdot,y)\phi(t,y)\dd y$ is a bounded linear operator on $C^0(\bar Q)$, and on $\mathcal O$ the denominator in \eqref{eq:Phi} is bounded away from zero on a neighborhood of each point; hence $\Phi$ is $C^1$ on $\mathcal O$ with locally Lipschitz derivative, and so is $\F$. Zeros of $\F$ are the mild solutions of \eqref{eq:mfgc} in $X$, and every classical solution is one of them. Writing $D\F(x^{n-1})(x^n-x^{n-1})=-\F(x^{n-1})$ componentwise shows that the Newton iterations for $\F$ are the mild formulation of \eqref{eq:newton}.

\emph{Step 2 (the control-law block).} By (v), $|D^2_{pZ}H|\le\lambda_0$. Since $m^*\ge0$ and $\eps>0$, for every $\zeta\in C^0(\bar Q;\mathbb R^d)$,
\[
|\partial_Z\Phi^*\zeta|\le\lambda_0\|\zeta\|_\infty\frac{\int Km^*}{\int Km^*+\eps}\le\lambda_0\|\zeta\|_\infty,
\]
so that $\|\partial_Z\Phi^*\|\le\lambda_0<1$ and $I-\partial_Z\Phi^*$ is invertible on $C^0(\bar Q;\mathbb R^d)$.

\emph{Step 3 ($D\F(x^*)$ is an isomorphism).} In block form with respect to $X=X_1\times X_2$, $X_1=C^{1,0}(\bar Q)\times C^0(\bar Q)$, $X_2=C^0(\bar Q;\mathbb R^d)$,
\[
D\F(x^*)=\begin{pmatrix} I-\mathcal K & -\mathcal B\\ -\mathcal C & I-\partial_Z\Phi^*\end{pmatrix},
\]
with
\[
\begin{aligned}
\mathcal K(v,\rho)&=\Big(\Tb\big[f'(m^*)\rho-H^*_p\cdot Dv\big],\ \Tf\big[\operatorname{div}\big(\rho H^*_p+m^*H^*_{pp}Dv\big)\big]\Big),\\
\mathcal B\zeta&=\Big(-\Tb\big[H^*_Z\cdot\zeta\big],\ \Tf\big[\operatorname{div}\big(m^*H^*_{pZ}\zeta\big)\big]\Big),\qquad
\mathcal C(v,\rho)=\partial_u\Phi^*v+\partial_m\Phi^*\rho .
\end{aligned}
\]
The arguments of $\Tb$ and $\Tf\circ\operatorname{div}$ depend continuously on $(v,\rho,\zeta)\in X$ in the uniform norm, so that $\mathcal K:X_1\to X_1$ and $\mathcal B:X_2\to X_1$ are compact by Lemma~\ref{lem:heat}, and $\mathcal C:X_1\to X_2$ is bounded. By Step 2, $D\F(x^*)$ is an isomorphism if and only if the Schur complement $\mathcal S=I-\mathcal K-\mathcal B(I-\partial_Z\Phi^*)^{-1}\mathcal C$ is an isomorphism of $X_1$. Since $\mathcal S$ is a compact perturbation of the identity, it is Fredholm of index $0$, and it is an isomorphism if and only if it is injective. If $\mathcal S(v,\rho)=0$, then $(v,\rho,\zeta)$ with $\zeta=(I-\partial_Z\Phi^*)^{-1}\mathcal C(v,\rho)$ belongs to $\ker D\F(x^*)$, i.e.\ it is a mild solution of \eqref{eq:lin} in $X$. By stability, $(v,\rho,\zeta)=0$. Note that the control-law block is not compact; it is the contraction (v) that allows its elimination.

\emph{Step 4 (conclusion).} $\F$ is $C^1$ near $x^*$ with locally Lipschitz derivative, $\F(x^*)=0$ and $D\F(x^*)$ is an isomorphism of $X$. By the local convergence theorem for Newton's method in Banach spaces (see \cite[Chapter~2]{Deu04} or \cite[Theorem 5.7]{BLS}), there is $\eta>0$ such that the Newton iterates starting in the ball of radius $\eta$ around $x^*$ are well defined and converge quadratically. Unique solvability of each step follows since $D\F(x)$ remains an isomorphism near $x^*$.
\end{proof}

Stable solutions exist when the interaction through the controls is weak. Write $H_\lambda=H_0(x,p)+\lambda\,h(x,p,Z)$, where $H_0$ and $h$ are of class $C^3$ with locally Lipschitz third derivatives, $D^2_{pp}H_0>0$, and $|D_ph(x,p,Z_1)-D_ph(x,p,Z_2)|\le|Z_1-Z_2|$; no convexity is required of $h$ (for the alignment model, $h=-Z\cdot p$). Then $H_\lambda$ satisfies (v) with $\lambda_0=|\lambda|$; the proof of Theorem~\ref{thm:main} uses only the smoothness in (ii) and (v), and the convexity enters below only at $\lambda=0$. As in Step 1 of that proof, $\F$ is of class $C^1$ in $(x,\lambda)$, with locally Lipschitz derivative.

\begin{corollary}\label{cor:small}
Assume (i), (iii), (iv), let $H=H_\lambda$ be as above, assume $\partial_mf\ge0$, and let $x_0^*$ be a classical solution of \eqref{eq:mfgc} with $\lambda=0$. Then there exist $\lambda^*\in(0,1)$ and a $C^1$ curve $\lambda\mapsto x^*_\lambda$, with $x^*_\lambda|_{\lambda=0}=x_0^*$ and $\|x^*_\lambda-x^*_0\|_X=O(\lambda)$, such that for every $\lambda\in(-\lambda^*,\lambda^*)$, $x^*_\lambda$ is a stable solution of \eqref{eq:mfgc}. Consequently, the conclusion of Theorem~\ref{thm:main} holds at $x^*_\lambda$ for all $\lambda\in(-\lambda^*,\lambda^*)$.
\end{corollary}

\begin{proof}
For $\lambda=0$ the system \eqref{eq:lin} is triangular: its first two equations are the linearized MFG system, whose only solution is $(0,0)$ when $\partial_mf\ge0$ and $D^2_{pp}H_0>0$ (energy identity, see \cite[Lemma 2.5]{CT}, using $m^*>0$ in $(0,T]\times\T$); the third equation then gives $\zeta=\partial_Z\Phi^*\zeta=0$, since $\partial_Z\Phi^*=0$ for $\lambda=0$. Hence $x_0^*$ is stable and $D_x\F(x^*_0;0)$ is an isomorphism by Step 3. Since $\F$ is $C^1$ in $(x,\lambda)$, the implicit function theorem gives the curve $x^*_\lambda$; by continuity, $D_x\F(x^*_\lambda;\lambda)$ remains an isomorphism for $\lambda$ small, and $x^*_\lambda$ is stable. The zero $x^*_\lambda\in X$ of $\F(\cdot;\lambda)$ is a mild solution of \eqref{eq:mfgc}; a bootstrap based on Lemma~\ref{lem:heat} and the parabolic Schauder estimates \cite[Chapter IV]{LSU} shows that it is classical, and $m^*_\lambda\ge0$ as in the proof of Theorem~\ref{thm:reduced}(a).
\end{proof}

\begin{remark}
Condition \eqref{eq:M} below ensures stability without requiring global uniqueness. In non-monotone problems, Theorem~\ref{thm:main} applies locally near each stable solution, including when several solutions coexist, as in \cite{BLS,AK21}.
\end{remark}

\subsection{Monotone interactions}

Corollary~\ref{cor:small} guarantees stability for sufficiently small $\lambda$, without providing an explicit threshold. We now show that, under an infinitesimal monotonicity condition of Lasry--Lions type for MFGC, a solution is stable without a further smallness condition on the interaction.

Differentiating $D_\alpha L(x,\alpha^*,Z)=-p$ and using the envelope theorem in \eqref{eq:legendreH} gives, at $(x,\alpha^*,Z)$,
\begin{equation}\label{eq:legendre}
D^2_{\alpha\alpha}L=(D^2_{pp}H)^{-1},\qquad D^2_{\alpha Z}L=(D^2_{pp}H)^{-1}D^2_{pZ}H,\qquad D_ZL=-D_ZH .
\end{equation}
The interaction depends on the joint law $\mu=(\mathrm{Id},\alpha^*)_\#m$ of states and controls through $Z=\mathcal Z[\mu]$; for \eqref{eq:Phi}, $\mathcal Z[\mu](x)=\int K(x,y)\,a\,\dd\mu(y,a)\,/\,(\int K(x,y)\dd\mu(y,a)+\eps)$. For a solution $(v,\rho,\zeta)$ of \eqref{eq:lin}, set
\[
\delta\alpha=-\big(H_{pp}^*Dv+H_{pZ}^*\zeta\big),
\]
the linearization of the optimal control, and let $\delta\mu(t)$ be the corresponding linearization of the joint law, acting on test functions $\varphi(x,a)$ by
\begin{equation}\label{eq:dmu}
\langle\delta\mu,\varphi\rangle=\int_{\T}\varphi(x,\alpha^*)\,\rho\dd x+\int_{\T}m^*\,D_a\varphi(x,\alpha^*)\cdot\delta\alpha\dd x .
\end{equation}
By the chain rule, the third equation of \eqref{eq:lin} is exactly $\zeta=D\mathcal Z[\mu^*](\delta\mu)$.

\begin{lemma}[Energy identity]\label{lem:energy}
Assume \eqref{eq:legendreH} with $L$ of class $C^3$, and let $(v,\rho,\zeta)\in X$ solve \eqref{eq:lin}. Then
\begin{equation}\label{eq:energy}
\int_Q m^*\,D^2_{\alpha\alpha}L\,\delta\alpha\cdot\delta\alpha\dd x\dd t+\int_Q f'(m^*)\rho^2\dd x\dd t+\int_0^T\big\langle\delta\mu(t),\,D_ZL(\cdot,\cdot,Z^*)\cdot\zeta(t)\big\rangle\dd t=0 .
\end{equation}
\end{lemma}
\begin{proof}
With $F$ and $g$ as in \eqref{eq:Fg}, the weak formulations of the two equations give $\langle\partial_t\rho,v\rangle=-\nu\int_{\T}D\rho\cdot Dv-\int_{\T}F\cdot Dv$ and $\langle\partial_tv,\rho\rangle=\nu\int_{\T}Dv\cdot D\rho-\int_{\T}g\rho$ for almost every $t$. Integrating their sum over $(0,T)$, and using the absolute continuity of $t\mapsto\int_{\T}\rho v$ and $\rho(0)=v(T)=0$, we obtain $\int_Q(F\cdot Dv+g\rho)=0$. The transport terms $\int\rho H_p^*\cdot Dv$ cancel, and this is
\[
\int_Q m^*\big(H^*_{pp}Dv+H^*_{pZ}\zeta\big)\cdot Dv-\int_Q (H_Z^*\cdot\zeta)\rho+\int_Q f'(m^*)\rho^2=0 .
\]
Writing $Dv=-(H^*_{pp})^{-1}(\delta\alpha+H^*_{pZ}\zeta)$ in the first integral and using \eqref{eq:legendre},
\[
\begin{aligned}
\int_Q m^*\big(H^*_{pp}Dv+H^*_{pZ}\zeta\big)\cdot Dv&=\int_Q m^*\,L_{\alpha\alpha}\delta\alpha\cdot\delta\alpha+\int_Q m^*\,\delta\alpha\cdot L_{\alpha Z}\zeta,\\
-\int_Q(H^*_Z\cdot\zeta)\rho&=\int_Q (L_Z\cdot\zeta)\rho .
\end{aligned}
\]
By \eqref{eq:dmu} with $\varphi(x,a)=D_ZL(x,a,Z^*(x))\cdot\zeta(x)$, for which $D_a\varphi=L_{\alpha Z}\zeta$, the last two terms are $\langle\delta\mu,D_ZL\cdot\zeta\rangle$.
\end{proof}

We assume the infinitesimal monotonicity condition: for all $t\in(0,T]$ and all $\delta\mu$ of the form \eqref{eq:dmu}, with $\rho\in C^0(\T)$, $\int_{\T}\rho=0$, and $\delta\alpha\in C^0(\T;\mathbb R^d)$ (the class of perturbations generated by solutions of \eqref{eq:lin}, since $\int\rho(t)=0$ by the divergence form of the second equation),
\begin{equation}\label{eq:M}
\big\langle\delta\mu,\ D_ZL(\cdot,\cdot,Z^*)\cdot D\mathcal Z[\mu^*](\delta\mu)\big\rangle\ge0. \tag{M}
\end{equation}
Condition \eqref{eq:M} is the second-order form of the monotonicity assumption of Kobeissi \cite{Kob20},
\[
\int\big(L(x,a,\mathcal Z[\mu^1])-L(x,a,\mathcal Z[\mu^2])\big)\dd(\mu^1-\mu^2)(x,a)\ge0\qquad\forall\,\mu^1,\mu^2,
\]
obtained by taking $\mu^1=(\mathrm{Id},\alpha^*+s\,\delta\alpha)_\#(m^*+s\rho)$, which is a probability measure for $|s|$ small when $m^*>0$, $\mu^2=\mu^*$, dividing by $s^2$ and letting $s\to0$: since $K$ and $L$ are smooth, $s\mapsto\mathcal Z[\mu^1]$ is twice differentiable, the integral equals $s^2\langle\delta\mu,D_ZL\cdot D\mathcal Z[\mu^*](\delta\mu)\rangle+o(s^2)$, and the global condition implies \eqref{eq:M}.

\begin{proposition}[Monotone case]\label{prop:monotone}
Let $x^*$ be a classical solution of \eqref{eq:mfgc}. Assume (i)--(v), \eqref{eq:legendreH} with $L$ of class $C^3$, $\partial_mf\ge0$, and \eqref{eq:M} evaluated at $x^*$. Then $x^*$ is stable. Consequently, Newton iterations initialized sufficiently close to $x^*$ converge quadratically.
\end{proposition}
\begin{proof}
By Lemma~\ref{lem:energy} and \eqref{eq:M}, $\int_Q m^*L_{\alpha\alpha}\delta\alpha\cdot\delta\alpha\le0$. Since $L_{\alpha\alpha}=(D^2_{pp}H)^{-1}>0$ by (ii) and $m^*>0$ in $(0,T]\times\T$ (strong maximum principle), $\delta\alpha=0$. The second equation of \eqref{eq:lin} reads $\partial_t\rho-\nu\Delta\rho-\operatorname{div}(\rho H^*_p)=-\operatorname{div}(m^*\delta\alpha)=0$ with $\rho(0)=0$, hence $\rho=0$. Then $\delta\mu=0$ and $\zeta=D\mathcal Z[\mu^*](\delta\mu)=0$, so $H^*_{pp}Dv=-\delta\alpha-H^*_{pZ}\zeta=0$ and $v=v(t)$; the first equation gives $-\partial_tv=0$ with $v(T)=0$, i.e.\ $v=0$.
\end{proof}
If the global monotonicity condition of \cite{Kob20} holds, the derivation above shows that \eqref{eq:M} holds at every classical solution, and hence every such solution is stable.

\begin{remark}[A normalized anti-alignment example]\label{rem:anti}
Consider the alignment Lagrangian \eqref{eq:align} with $K(x,y)\equiv1$
on the unit-volume torus. For probability densities,
\[
\mathcal Z[\mu]=\frac{\int a\,\mathrm d\mu(x,a)}{1+\varepsilon}.
\]
For a mass-preserving perturbation \eqref{eq:dmu}, set
\[
A=\int_{\T}(\alpha^*\rho+m^*\delta\alpha)\,\dd x,
\qquad \int_{\T}\rho\,\dd x=0.
\]
Then $D\mathcal Z[\mu^*](\delta\mu)=A/(1+\varepsilon)$ and
$D_ZL(x,\alpha,Z)=-\lambda(\alpha-\lambda Z)$, so
\begin{equation}\label{eq:anti-normalized}
\big\langle\delta\mu,\,
D_ZL(\cdot,\cdot,Z^*)\cdot D\mathcal Z[\mu^*](\delta\mu)
\big\rangle
=-\frac{\lambda}{1+\varepsilon}|A|^2.
\end{equation}
Thus condition \eqref{eq:M} holds for $-1<\lambda\le0$,
providing an interacting example within the normalized model
\eqref{eq:Phi}. For $\lambda>0$ the sign is reversed. For non-constant $K$, $D\mathcal Z[\mu^*](\delta\mu)$ depends on $x$, and the sign in \eqref{eq:M} is not determined by the sign of $\lambda$ alone; we do not pursue this here.
\end{remark}

\subsection{Newton iterations on reduced unknowns}\label{sec:reduced}

The Newton iterations \eqref{eq:newton} update the density by the solution $m^n$ of a linear FP equation with the source $\operatorname{div}\big(m^{n-1}[H_{pp}^{n-1}(Du^n-Du^{n-1})+H^{n-1}_{pZ}(Z^n-Z^{n-1})]\big)$, which has no prescribed sign. The iterate $m^n$ need not be nonnegative, so that the denominator in \eqref{eq:Phi} may vanish along the iterations. We now consider Newton iterations on the reduced unknowns $(g,Z)$, where $g$ stands for the coupling term $f(\cdot,m)$, in which the value function and the density are always exact solutions of their own equations. This is the continuous counterpart of the strategy of \cite{AK21}, and an instance of nonlinear elimination in Newton methods \cite{LRW96}.

For $(g,Z)\in \mathcal Y=C^0(\bar Q)\times C^0(\bar Q;\mathbb R^d)$, let $u[g,Z]$ solve
\begin{equation}\label{eq:hjbred}
-\partial_tu-\nu\Delta u+H(x,Du,Z)=g\quad\text{in }Q,\qquad u(T)=G,
\end{equation}
let $m[g,Z]$ solve
\begin{equation}\label{eq:fpred}
\partial_tm-\nu\Delta m-\operatorname{div}\big(m\,D_pH(x,Du[g,Z],Z)\big)=0\quad\text{in }Q,\qquad m(0)=m_0,
\end{equation}
and define
\begin{equation}\label{eq:Psi}
\Psi(g,Z)=\Big(f\big(\cdot,m[g,Z]\big),\ \Phi\big(u[g,Z],m[g,Z],Z\big)\Big),\qquad \mathcal R(g,Z)=(g,Z)-\Psi(g,Z).
\end{equation}
Then $(u,m,Z)$ solves \eqref{eq:mfgc} if and only if $\mathcal R(g,Z)=0$ with $g=f(\cdot,m)$, $u=u[g,Z]$, $m=m[g,Z]$. The Newton iterations on the reduced unknowns are
\begin{equation}\label{eq:newtonred}
D\mathcal R(g^{n-1},Z^{n-1})\,(g^n-g^{n-1},Z^n-Z^{n-1})=-\mathcal R(g^{n-1},Z^{n-1}).
\end{equation}
The derivative of $\Psi$ in a direction $(\delta g,\delta Z)$ is obtained by solving one linearized HJB equation backward and one linearized FP equation forward,
\begin{equation}\label{eq:tangent}
\begin{cases}
-\partial_t\delta u-\nu\Delta\delta u+H_p\cdot D\delta u=\delta g-H_Z\cdot\delta Z, & \delta u(T)=0,\\[2pt]
\partial_t\delta m-\nu\Delta\delta m-\operatorname{div}(\delta m\,H_p)=\operatorname{div}\big(m\,[H_{pp}D\delta u+H_{pZ}\delta Z]\big), & \delta m(0)=0,
\end{cases}
\end{equation}
and then $D\Psi(\delta g,\delta Z)=(f'(m)\delta m,\ \partial_u\Phi\,\delta u+\partial_m\Phi\,\delta m+\partial_Z\Phi\,\delta Z)$. The equations \eqref{eq:tangent} are the linearized equations of \eqref{eq:newton}; the difference with \eqref{eq:newton} is that the state $(u,m)$ is not updated by the solution of the linear system, but recomputed by \eqref{eq:hjbred}--\eqref{eq:fpred} at the new reduced unknowns.

\begin{theorem}\label{thm:reduced}
Assume (i)--(v), let $x^*=(u^*,m^*,Z^*)$ be a stable classical solution of \eqref{eq:mfgc}, and let $y^*=(f(\cdot,m^*),Z^*)$. Then:
\begin{enumerate}[(a)]
\item there is a neighborhood $\mathcal U$ of $y^*$ in $\mathcal Y$ such that $u[g,Z]\in C^{1,0}(\bar Q)$ and $m[g,Z]\in C^0(\bar Q)$ are well defined for $(g,Z)\in\mathcal U$, depend in a $C^1$ way on $(g,Z)$ with locally Lipschitz derivative, and $m[g,Z]\ge0$, $\int_{\T}m[g,Z](t)\dd x=1$ for all $t$;
\item $D\mathcal R(y^*)$ is an isomorphism of $\mathcal Y$;
\item there exist $\eta>0$ and $C>0$ such that, if $\|y^0-y^*\|_\mathcal Y\le\eta$, the iterations \eqref{eq:newtonred} are well defined and $\|y^n-y^*\|_\mathcal Y\le C\|y^{n-1}-y^*\|_\mathcal Y^2$ for all $n\ge1$.
\end{enumerate}
\end{theorem}
\begin{proof}
Throughout, $u[g,Z]$ and $m[g,Z]$ are mild solutions, i.e.\ $u=\bar u+\Tb[g-H(\cdot,Du,Z)]$ and $m=\bar m+\Tf[\operatorname{div}(mD_pH(\cdot,Du,Z))]$.

(a) Consider $E_1(u;g,Z)=u-\bar u-\Tb[g-H(\cdot,Du,Z)]$, a map of class $C^2$ from $C^{1,0}(\bar Q)\times \mathcal Y$ to $C^{1,0}(\bar Q)$, by Lemma~\ref{lem:heat} and since the superposition operator associated with $H\in C^3$ is of class $C^2$ on spaces of continuous functions \cite[Chapter 1]{AP93}. Its derivative in $u$ at $(u^*;y^*)$ is $v\mapsto v+\Tb[H^*_p\cdot Dv]$, a compact perturbation of the identity, which is injective by the remark following Lemma~\ref{lem:heat}; hence it is an isomorphism. By the implicit function theorem, $u[g,Z]$ is defined near $y^*$, with values in $C^{1,0}(\bar Q)$, and is of class $C^2$, so that its derivative is locally Lipschitz. Next, for a drift $b\in C^0(\bar Q;\mathbb R^d)$, let $\mathcal K_b\rho=\Tf[\operatorname{div}(\rho b)]$. By Lemma~\ref{lem:heat}, $\mathcal K_b$ is a compact operator on $C^0(\bar Q)$, and $I-\mathcal K_b$ is injective by the same remark; hence $I-\mathcal K_b$ is an isomorphism of $C^0(\bar Q)$, and $m=(I-\mathcal K_b)^{-1}\bar m$ is the unique solution of $m=\bar m+\Tf[\operatorname{div}(mb)]$. Since $b\mapsto\mathcal K_b$ is linear and continuous with values in the bounded operators on $C^0(\bar Q)$, by Lemma~\ref{lem:heat}(a), and the inversion of operators is of class $C^\infty$, $m$ depends in a $C^\infty$ way on $b$. Since $(g,Z)\mapsto D_pH(\cdot,Du[g,Z],Z)\in C^0(\bar Q;\mathbb R^d)$ is $C^1$ with locally Lipschitz derivative, so is $m[g,Z]$. Finally, $\int_{\T}m[g,Z](t)=\int_{\T}m_0=1$ by Lemma~\ref{lem:heat}(c), and $m[g,Z]\ge0$: approximating the drift $b$ uniformly by smooth drifts $b_k$, the corresponding solutions $m_k$ are classical and nonnegative by the maximum principle, and $\rho_k=m_k-m$ satisfies $\rho_k=\Tf[\operatorname{div}(\rho_kb)]+\Tf[\operatorname{div}(m_k(b_k-b))]$. Since the drifts $b_k$ are bounded uniformly in $k$, Lemma~\ref{lem:heat}(a) and the Gronwall lemma applied to $m_k=\bar m+\Tf[\operatorname{div}(m_kb_k)]$ give a bound $\|m_k\|_\infty\le C_0$ independent of $k$. Hence, by Lemma~\ref{lem:heat}(a),
\[
\|\rho_k(t)\|_\infty\le C\|b\|_\infty\int_0^t\frac{\|\rho_k(s)\|_\infty}{(t-s)^{1/2}}\dd s+\eps_k,\qquad \eps_k=2CT^{1/2}C_0\|b_k-b\|_\infty\to0,
\]
and the Gronwall lemma gives $\|\rho_k\|_\infty\le C'\eps_k\to0$. Thus $m=\lim m_k\ge0$.

(b) Consider the extended map $\mathcal E(u,m,g,Z)=(E_1(u;g,Z),\,E_2(u,m;Z),\,g-f(\cdot,m),\,Z-\Phi(u,m,Z))$, with $E_2(u,m;Z)=m-\bar m-\Tf[\operatorname{div}(mD_pH(\cdot,Du,Z))]$. Eliminating $g$ by the third component shows that $D\mathcal E$ at $(u^*,m^*,y^*)$ is an isomorphism if and only if $D\F(x^*)$ is, which holds by Step 3 of the proof of Theorem~\ref{thm:main}. On the other hand, write $D\mathcal E=\begin{pmatrix}\mathcal E_{11}&\mathcal E_{12}\\ \mathcal E_{21}&\mathcal E_{22}\end{pmatrix}$ with respect to the splitting into $(u,m)$ and $(g,Z)$ for the unknowns, and into $(E_1,E_2)$ and the last two components for the equations. The block $\mathcal E_{11}$ is block triangular (the derivative of $E_1$ in $m$ vanishes) with invertible diagonal blocks, by (a), hence invertible. Differentiating $E_1(u[y];y)=0$ and $E_2(u[y],m[y];y)=0$ gives $D(u,m)[y^*]=-\mathcal E_{11}^{-1}\mathcal E_{12}$, so that $D\mathcal R(y^*)=\mathcal E_{22}-\mathcal E_{21}\mathcal E_{11}^{-1}\mathcal E_{12}$ is the Schur complement of $\mathcal E_{11}$. A block operator with invertible $\mathcal E_{11}$ is an isomorphism if and only if its Schur complement is; hence $D\mathcal R(y^*)$ is an isomorphism.

(c) By (a), $\mathcal R$ is $C^1$ near $y^*$ with locally Lipschitz derivative, $\mathcal R(y^*)=0$, and $D\mathcal R(y^*)$ is an isomorphism by (b). We conclude by the local convergence theorem for Newton's method, as in Step 4 of the proof of Theorem~\ref{thm:main}.
\end{proof}

\section{A semi-Lagrangian scheme for MFGC}\label{sec:sl}

For $N,N_t\in\mathbb N$, let $h=1/N$ be the spatial mesh size, $\mathcal G_h=\{x_i\}=\{hj:\ j\in\{0,\dots,N-1\}^d\}$ the uniform grid on $\T=\mathbb R^d/\mathbb Z^d$, which has $N^d$ nodes, $\Delta t=T/N_t$ the time step and $t_k=k\Delta t$. We denote by $\I[\phi]$ the $\mathbb Q_1$ interpolant of a grid function $\phi$, by $D_h$ the centered discrete gradient and by $\operatorname{div}_h=-D_h^\top$ the centered discrete divergence, and by $K_h$ the discrete kernel, $(K_h\phi)_i=h^d\sum_jK(x_i,x_j)\phi_j$. Grid functions at time $t_k$ are denoted by $U^k$, $M^k$, $Z^k$. We present the scheme on the torus; the modifications for a square with reflecting boundary are given at the end of this section.

\subsection{The HJB equation}
Following \cite{CS15,CZ}, the value function is approximated along the characteristics of the optimal control. For a grid function $U^k$ and $Z^k$, let
\[
\alpha^k_i=-D_pH\big(x_i,(D_hU^k)_i,Z^k_i\big),\qquad y^{k,\ell}_i=x_i+\Delta t\,\alpha^k_i+\sqrt{2d\,\nu\Delta t}\,e^\ell,\quad \ell=1,\dots,2d,
\]
where $e^\ell\in\{\pm e_1,\dots,\pm e_d\}$, and define
\begin{equation}\label{eq:slhjb}
U^k_i=\frac1{2d}\sum_{\ell=1}^{2d}\I[U^{k+1}]\big(y^{k,\ell}_i\big)+\Delta t\Big(L\big(x_i,\alpha^k_i,Z^k_i\big)+f(x_i,M^k_i)\Big),\qquad U^{N_t}_i=G(x_i).
\end{equation}
The displacements $\xi^\ell=\sqrt{2d\nu\Delta t}\,e^\ell$ satisfy $\frac1{2d}\sum_\ell\xi^\ell=0$ and $\frac1{2d}\sum_\ell\xi^\ell(\xi^\ell)^\top=2\nu\Delta t\,I$, so that, for smooth $\phi$, $\frac1{2d}\sum_\ell\phi(x+\xi^\ell)=\phi(x)+\nu\Delta t\,\Delta\phi(x)+O(\Delta t^2)$: the average reproduces the diffusion term of $\dd X=\alpha\dd t+\sqrt{2\nu}\dd W$, as in \cite{CS15,CZ}. For the alignment model \eqref{eq:align}, $\alpha^k_i=\lambda Z^k_i-(D_hU^k)_i$ and $L(x_i,\alpha^k_i,Z^k_i)=\frac12|(D_hU^k)_i|^2$. Scheme \eqref{eq:slhjb} is the SL scheme of \cite{CZ} at its fixed point, i.e.\ when the drift computed from the previous Newton iterate coincides with the current one, extended to Hamiltonians depending on $Z$. It is implicit only through the discrete gradient defining the feet of the characteristics. Unlike the scheme of \cite{CS15}, in which the control is obtained by a minimization over $\alpha$, the control in \eqref{eq:slhjb} is given by the centered gradient $D_hU^k$; the scheme is therefore not monotone in general, and the convergence results of \cite{CS15} do not apply to it.

\subsection{The FP equation}
Writing \eqref{eq:slhjb} as $U^k=B^k[U^k,Z^k]\,U^{k+1}+\Delta t(\dots)$, where $B^k$ is the matrix of interpolation weights,
\[
\big(B^k\phi\big)_i=\frac1{2d}\sum_{\ell=1}^{2d}\I[\phi]\big(y^{k,\ell}_i\big),
\]
the discrete FP equation is defined with its transpose,
\begin{equation}\label{eq:slfp}
M^{k+1}=\big(B^k\big)^\top M^k,\qquad M^0_i=\frac{m_0(x_i)}{h^d\sum_jm_0(x_j)}.
\end{equation}
The normalization of $M^0$ gives $h^d\sum_iM^0_i=1$. In practice $(B^k)^\top M^k$ is computed without assembling $B^k$: each node distributes its mass $M^k_i$ to the vertices of the cells containing the feet $y^{k,\ell}_i$, with the interpolation weights. Here $(B^k)^\top$ is the adjoint of the interpolation operator at the computed characteristics, not of the derivative of $U^k\mapsto B^k[U^k,Z^k]U^{k+1}$, which would also contain the derivatives of the weights with respect to the feet; this is the construction of \cite{CS15,CZ}, whereas in \cite{AC10,AK21} the discrete FP operator is the adjoint of the linearized discrete HJB operator.

\subsection{The average control}
The relation \eqref{eq:Phi} is discretized at each time level by
\begin{equation}\label{eq:slZ}
Z^k\,\big(K_hM^k+\eps\big)=K_h\big[\alpha^k\,M^k\big],
\end{equation}
where products are taken componentwise. For the alignment model, \eqref{eq:slZ} is linear in $Z^k$:
\begin{equation}\label{eq:slZlin}
\big(I-\lambda S^k\big)Z^k=-\big(K_hM^k+\eps\big)^{-1}K_h\big[M^k D_hU^k\big],\qquad S^k\zeta=\big(K_hM^k+\eps\big)^{-1}K_h\big[M^k\zeta\big].
\end{equation}

Equations \eqref{eq:slhjb}--\eqref{eq:slZ} define the discrete MFGC system. For a smooth solution, the formal local truncation error of the SL update for the HJB equation is $O(\Delta t+h^2/\Delta t)$ \cite{CS15,CZ}. The convergence results for classical MFG in \cite{CS15}, including the requirement $h^2/\Delta t\to0$, do not establish convergence of the present coupled MFGC scheme.

\begin{lemma}\label{lem:sl}
Let $M^0\ge0$ with $h^d\sum_iM^0_i=1$. For any $U^k$ and $Z^k$, the solution of \eqref{eq:slfp} satisfies $M^k\ge0$ and $h^d\sum_iM^k_i=1$ for all $k$.
\end{lemma}
\begin{proof}
The entries of $B^k$ are nonnegative, since the $\mathbb Q_1$ interpolation weights are, and each row of $B^k$ sums to one, since the interpolation reproduces constants. Hence $(B^k)^\top$ maps nonnegative vectors to nonnegative vectors and preserves the sum of the entries.
\end{proof}

\begin{lemma}\label{lem:Zdiscrete}
Let $M^k\ge0$ and assume that the matrix of $K_h$ is symmetric with nonnegative entries. Let $\mathcal P=\{i:M^k_i>0\}$. Then the spectrum of $S^k$ in \eqref{eq:slZlin} is real and contained in $[-s,s]$, where
\[
s=\max_i\frac{(K_hM^k)_i}{(K_hM^k)_i+\eps}<1,
\]
and in $[0,s]$ if $K_h$ is positive semidefinite. Consequently, for $|\lambda|<1$, \eqref{eq:slZlin} has a unique solution. Moreover, if $0<|\lambda|<1$, the Chebyshev iteration applied to the full system \eqref{eq:slZlin} with the interval $J=[1-|\lambda|,1+|\lambda|]$, which contains the spectrum of $I-\lambda S^k$, converges: with $\varepsilon_n=2\big((\sqrt\kappa-1)/(\sqrt\kappa+1)\big)^n$, $\kappa=(1+|\lambda|)/(1-|\lambda|)$, and $R=\operatorname{diag}\big((K_hM^k+\eps)^{1/2}(M^k)^{1/2}\big)_{\mathcal P}$, its errors $e_n$ satisfy
\[
\begin{gathered}\|R\,(e_n)_{\mathcal P}\|_2\le\varepsilon_n\|R\,(e_0)_{\mathcal P}\|_2,\\ \|(e_n)_{\mathcal P^c}\|_2\le\varepsilon_n\Big(\|(e_0)_{\mathcal P^c}\|_2+n^2\,\big\|S_{\mathcal P^c\mathcal P}R^{-1}\big\|_2\,\|R\,(e_0)_{\mathcal P}\|_2\Big).\end{gathered}
\]
\end{lemma}
\begin{proof}
We write $S^k=D_K^{-1}K_hD_M$, with $D_K=\operatorname{diag}(K_hM^k+\eps)$ and $D_M=\operatorname{diag}(M^k)$. Since the columns of $S^k$ indexed by $i\notin\mathcal P$ vanish, ordering the nodes as $(\mathcal P,\mathcal P^c)$ gives
\[
S^k=\begin{pmatrix}S_{\mathcal P\mathcal P}&0\\ S_{\mathcal P^c\mathcal P}&0\end{pmatrix},
\]
so that the spectrum of $S^k$ is that of $S_{\mathcal P\mathcal P}$ together with $0$. On $\mathcal P$ the diagonal matrices are invertible and, with $R=(D_K^{1/2}D_M^{1/2})_{\mathcal P}$,
\[
R\,S_{\mathcal P\mathcal P}\,R^{-1}=\big(D_K^{-1/2}D_M^{1/2}K_hD_M^{1/2}D_K^{-1/2}\big)_{\mathcal P\mathcal P},
\]
which is symmetric, and positive semidefinite if $K_h$ is. Hence the spectrum of $S^k$ is real; since $S^k$ has nonnegative entries, it is bounded by $\|S^k\|_\infty=s$. Since $|\lambda|s<1$, $I-\lambda S^k$ is invertible. The implementation applies the Chebyshev recurrence to the full system, so that $e_n=p_n(A)e_0$, where $A=I-\lambda S^k$ and $p_n$ is the Chebyshev polynomial of degree $n$ for $J$ normalized by $p_n(0)=1$, which satisfies $\max_J|p_n|\le\varepsilon_n$ \cite[Section 12.1]{Saa03}. In the ordering $(\mathcal P,\mathcal P^c)$, $A=\begin{pmatrix}A_1&0\\ A_2&I\end{pmatrix}$ with $A_1=I-\lambda S_{\mathcal P\mathcal P}$ and $A_2=-\lambda S_{\mathcal P^c\mathcal P}$, and by induction on the powers of $A$,
\[
p_n(A)=\begin{pmatrix}p_n(A_1)&0\\ A_2\,q_n(A_1)&p_n(1)I\end{pmatrix},\qquad q_n(t)=\frac{p_n(t)-p_n(1)}{t-1},
\]
where $q_n$ is a polynomial. Since $RA_1R^{-1}$ is symmetric with spectrum in~$J$, we have $\|Rp_n(A_1)R^{-1}\|_2\le\max_J|p_n|$ and \[\|Rq_n(A_1)R^{-1}\|_2\le\max_J|q_n|\le\max_J|p_n'|\le\frac{2n^2}{|J|}\max_J|p_n|,\] by the mean value theorem ($1\in J$) and Markov's inequality, with $|J|=2|\lambda|$, and $\|A_2\,\cdot\|=|\lambda|\|S_{\mathcal P^c\mathcal P}\,\cdot\|$; moreover $|p_n(1)|\le\max_J|p_n|$. The two bounds follow from \[(e_n)_{\mathcal P}=p_n(A_1)(e_0)_{\mathcal P},\qquad (e_n)_{\mathcal P^c}=A_2q_n(A_1)(e_0)_{\mathcal P}+p_n(1)(e_0)_{\mathcal P^c}.\]
\end{proof}

\subsection{Reflecting boundary conditions}
In Test 4 of Section 6, the domain is the square $\Omega=(-\frac12,\frac12)^d$ with homogeneous Neumann boundary conditions for $u$ and no-flux boundary conditions for $m$, which correspond to agents driven by a reflected Brownian motion. The grid is $\mathcal G_h=\{-\frac12+jh\}^d$, $j=0,\dots,N$, $h=1/N$, including the boundary nodes. The feet of the characteristics are reflected into $\bar\Omega$ coordinatewise by the map
\[
\pi(y)_c=-\tfrac12+\min\big(z_c,\,2-z_c\big),\qquad z_c=\big(y_c+\tfrac12\big)\bmod 2,\qquad c=1,\dots,d,
\]
which accounts for any number of reflections and for corners, and \eqref{eq:slhjb} is applied with $\I[U^{k+1}](\pi(y^{k,\ell}_i))$. The interpolation weights at the reflected feet are nonnegative and sum to one, so that $B^k$ has the same properties as on the torus and Lemma~\ref{lem:sl} holds for \eqref{eq:slfp} on $\mathcal G_h$: the mass that would leave $\Omega$ is reflected back, which is the discrete counterpart of the no-flux condition. The discrete gradient $D_h$ is centered in the interior and its normal component is set to zero at the boundary nodes, consistently with the Neumann condition, and $(K_h\phi)_i=h^d\sum_{x_j\in\mathcal G_h}K(x_i,x_j)\phi_j$. For the symmetric kernel $K(x,y)=K(y,x)$ used in Test 4, this matrix is symmetric and has nonnegative entries. Lemma~\ref{lem:Zdiscrete} therefore applies (with the interval $[-s,s]$ for a kernel which is not positive semidefinite), whereas Proposition~\ref{prop:bgs} below is proved on the torus only, where $\operatorname{div}_h=-D_h^\top$. The theoretical results of Section~\ref{sec:theory} are stated on the torus; we do not analyze the bounded domain.

\section{Newton solvers for the semi-Lagrangian scheme}\label{sec:solvers}

\subsection{Newton--SL iterations and block Gauss--Seidel}\label{sec:bgs}

We discretize the linearized HJB and FP equations of \eqref{eq:newton} by the SL scheme, as in \cite{CZ}. For the control equation, we apply Newton linearization to the denominator-cleared relation $Z(K*m+\eps)=K*[m\,\alpha^*]$. Its residual is the quotient residual in \eqref{eq:newton} multiplied by $K*m+\eps>0$; the two residuals have the same zeros, and at a zero their derivatives differ by an invertible factor, so the local argument of Theorem~\ref{thm:main} applies. For the alignment model, let $(U_o,M_o,Z_o)$ be the previous Newton iterate and $B^k$ the interpolation matrix associated with $q^k=D_hU_o^k-\lambda Z_o^k$. The new iterate $(U,M,Z)$ solves
{\small
\begin{equation}\label{eq:nsl}
\begin{cases}
U^k=B^kU^{k+1}+\Delta t\,\lambda\,D_hU_o^k\cdot Z^k+\Delta t\,W^kM^k+\Delta t\,b_U^k, & U^{N_t}=G,\\[2pt]
M^{k+1}=(B^k)^\top M^k+\Delta t\,\mathcal N^{k+1}U^{k+1}-\Delta t\,\lambda\operatorname{div}_h\big(M_o^{k+1}Z^{k+1}\big)+\Delta t\,b_M^{k+1}, & M^0=m_0,\\[2pt]
Z^k(K_hM_o^k+\eps)-\lambda K_h[M_o^kZ^k]+\mathcal A_m^kM^k+\mathcal A_u^kU^k=b_Z^k,
\end{cases}
\end{equation}}
where
\[
\begin{aligned}
W^k&=\operatorname{diag}(f'(M_o^k)),&
\mathcal N^kV&=\operatorname{div}_h(M_o^kD_hV),\\
\mathcal A_m^kM&=Z_o^kK_hM-K_h[(\lambda Z_o^k-D_hU_o^k)M],&
\mathcal A_u^kU&=K_h[M_o^kD_hU].
\end{aligned}
\]
The right-hand sides are explicit functions of the previous iterate. With $p_o^k=D_hU_o^k$ and $q_o^k=p_o^k-\lambda Z_o^k$, they are
\begin{equation}\label{eq:nsl-forcing}
\begin{aligned}
b_U^k&=f(\cdot,M_o^k)-f'(\cdot,M_o^k)M_o^k
                   +\tfrac12|p_o^k|^2-\lambda Z_o^k\cdot p_o^k,\\
b_M^k&=-\operatorname{div}_h(M_o^k q_o^k),\\
b_Z^k&=Z_o^k(K_hM_o^k)-\lambda K_h[M_o^k Z_o^k]
                         +K_h[M_o^k p_o^k].
\end{aligned}
\end{equation}
All products are taken pointwise on the grid. For $\lambda=0$, the first two equations coincide with the Newton--SL scheme of \cite{CZ}.

We solve \eqref{eq:nsl} by block Gauss--Seidel (BGS): a backward sweep for $U$, a forward sweep for $M$ and a solution of the control equation for $Z$. All operations are matrix-free. For $\lambda=0$, write the first two equations as
\[
\mathbf A U=\Delta t\,\mathbf W M+\mathbf b_U,\qquad
\mathbf A^\top M=\Delta t\,\mathbf N U+\mathbf b_M,
\]
where $\mathbf A$ has diagonal blocks $I$ and first superdiagonal blocks $-B^k$, and
\[
\mathbf W=\operatorname{diag}(W^0,\ldots,W^{N_t-1},0),\qquad \mathbf N=\operatorname{diag}(0,\mathcal N^1,\ldots,\mathcal N^{N_t}).
\]
More precisely,
\[
\mathbf b_U=(\Delta t\,b_U^0,\ldots,
\Delta t\,b_U^{N_t-1},G),\qquad
\mathbf b_M=(m_0,\Delta t\,b_M^1,\ldots,
\Delta t\,b_M^{N_t}).
\]

\begin{proposition}\label{prop:bgs}
Assume $\lambda=0$, $M_o\ge0$ and $f'\ge0$. Set
\[
Q=-\Delta t^2\,\mathbf N,\qquad
P=\mathbf A^{-\top}Q\mathbf A^{-1}\mathbf W,\qquad
\mu=\rho(P).
\]
The spectrum of the BGS error operator is real and contained in $[-\mu,0]$. Consequently:
\begin{enumerate}[(a)]
\item BGS converges for every right-hand side and initial guess if and only if $\mu<1$;
\item the relaxed iteration $V_{\ell+1}=V_\ell+\omega(\mathcal S(V_\ell)-V_\ell)$, with $\mathcal S$ denoting a BGS sweep, converges for
\[
0<\omega<\frac{2}{1+\mu}.
\]
The optimal choice $\omega^*=2/(2+\mu)$ yields the asymptotic convergence factor $\mu/(2+\mu)$.
\end{enumerate}
\end{proposition}
\begin{proof}
The BGS errors satisfy
\[
e_U^{\ell+1}=\Delta t\,\mathbf A^{-1}\mathbf W e_M^\ell,
\qquad e_M^{\ell+1}=-P e_M^\ell.
\]
Since $\operatorname{div}_h=-D_h^\top$ and $M_o\ge0$, the matrix $Q=\Delta t^2D_h^\top\operatorname{diag}(M_o)D_h$ is positive semidefinite. The nonzero eigenvalues of $P$ are those of
\[
\mathbf W^{1/2}\mathbf A^{-\top}Q\mathbf A^{-1}\mathbf W^{1/2},
\]
which is symmetric positive semidefinite. Thus $\sigma(P)\subset[0,\mu]$ and $\sigma(\mathcal S)\subset[-\mu,0]$. The relaxed eigenvalues are $1-\omega(1+\mu_j)$, $\mu_j\in[0,\mu]$, including $1-\omega$ for zero eigenvalues. Minimizing their maximum modulus gives the stated range and optimal value of $\omega$.
\end{proof}

For $\lambda\ne0$, the spectrum need not lie on the nonpositive real axis (Section~\ref{sec:test1}). We also use GMRES on the BGS-preconditioned system, with one sweep per matrix--vector product \cite{SS86}. The outer full Newton iterations follow \cite[Algorithm~2]{CZ}, with an Armijo line search when needed. Thus the full method linearizes and then discretizes, whereas the reduced method of the next subsection differentiates its discrete residual approximately. As noted in \cite{CZ} for MFG, the SL discretization of the Newton iterations is not the Newton method of the discrete system, since the derivatives of the interpolation weights with respect to the feet are not included; quadratic convergence at the discrete level is therefore not expected, and in Test~1 we observe linear convergence.

\subsection{Reduced Newton--SL iterations}\label{sec:ne}

Let $Y=(\Gamma,Z)$, where $\Gamma=(\Gamma^k)_{k<N_t}$ and $Z=(Z^k)_{k\le N_t}$. Given $Y$, compute $U=U[Y]$ backward from
\begin{equation}\label{eq:reducedsl-hjb}
\begin{aligned}
U_i^k&=\frac1{2d}\sum_{\ell=1}^{2d}\I[U^{k+1}]\!\left(x_i+\Delta t\,\alpha_i^k+\sqrt{2d\nu\Delta t}\,e^\ell\right)
+\Delta t\left(L(x_i,\alpha_i^k,Z_i^k)+\Gamma_i^k\right),\\
\alpha_i^k&=-D_pH\big(x_i,(D_hU^k)_i,Z_i^k\big),\qquad U^{N_t}=G,
\end{aligned}
\end{equation}
by a fixed-point iteration initialized with $U^{k,(0)}=U^{k+1}$ at each time level. We do not analyze the solvability of \eqref{eq:reducedsl-hjb} or the existence of discrete solutions; the Newton driver rejects every iterate for which a fixed-point iteration or the Chebyshev iteration of Lemma~\ref{lem:Zdiscrete} does not reach its tolerance. Then compute $M=M[Y]$ and $\widetilde Z=\widetilde Z[Y]$ from
\begin{equation}\label{eq:reducedsl-forward}
\begin{gathered}
M^{k+1}=(B^k)^\top M^k,\qquad M^0\text{ as in \eqref{eq:slfp}},\\
\big(I-\lambda S^k(M^k)\big)\widetilde Z^k
=-\big(K_hM^k+\eps\big)^{-1}K_h[M^kD_hU^k],
\end{gathered}
\end{equation}
where the control equation is solved as in Lemma~\ref{lem:Zdiscrete}. Define
\begin{equation}\label{eq:reducedsl-residual}
\Psi_h(Y)=\big((f(\cdot,M^k))_{k<N_t},\widetilde Z\big),
\qquad
\mathcal R_h(Y)=Y-\Psi_h(Y).
\end{equation}
The zeros of $\mathcal R_h$ solve \eqref{eq:slhjb}--\eqref{eq:slZ}. By Lemma~\ref{lem:sl} and the normalization of $M^0$, the reconstructed density satisfies $M^k\ge0$ and $h^d\sum_iM_i^k=1$ at every iteration.

\begin{remark}[Relation with \eqref{eq:Psi}]\label{rem:elim}
The implemented map solves the control equation for $\widetilde Z$, rather than evaluating $\Phi$ at the input $Z$. For the alignment model, let $A(M)=I-\lambda S(M)$ and $\ell(U,M)=(K_hM+\eps)^{-1}K_h[M D_hU]$. Then
\[
r_Z=A(M)Z+\ell(U,M),\qquad Z-\widetilde Z=A(M)^{-1}r_Z.
\]
Thus the two formulations have the same zeros. At a zero, their Jacobians differ by multiplication with $\operatorname{diag}(I,A(M)^{-1})$. Similarly, for the corresponding continuous map $\widetilde{\mathcal R}$,
\[
D\widetilde{\mathcal R}(y^*)=
\operatorname{diag}\big(I,(I-\partial_Z\Phi^*)^{-1}\big)D\mathcal R(y^*).
\]
Hence Theorem~\ref{thm:reduced} also applies locally to this formulation; away from a zero, the Newton iterates need not coincide.
\end{remark}

We solve $\mathcal R_h(Y)=0$ by Jacobian-free Newton--GMRES:
\begin{equation}\label{eq:reducedsl-newton}
D\mathcal R_h(Y^n)d^n=-\mathcal R_h(Y^n),
\qquad Y^{n+1}=Y^n+a_nd^n,
\end{equation}
with finite-difference Jacobian--vector products \cite{KK04} (the product with $v=0$ is zero),
\begin{equation}\label{eq:reducedsl-jv}
D\mathcal R_h(Y)v\approx
\frac{\mathcal R_h(Y+\sigma v)-\mathcal R_h(Y)}{\sigma},
\qquad \sigma=\sigma_0\frac{1+\|Y\|_2}{\|v\|_2},\quad v\neq0,
\end{equation}
with $\sigma_0>0$. On a region where $D\mathcal R_h$ is Lipschitz continuous with constant $L$, if the residual is evaluated with an error at most $\delta$, the error of \eqref{eq:reducedsl-jv} is at most $\frac L2\sigma\|v\|_2^2+2\delta/\sigma$. This estimate does not apply across the kinks of the interpolation weights.
We first take $a_n=1$. If this undamped phase produces a nonfinite residual, a failed inner solve, or $\Theta_h(Y^n)>\gamma\,\Theta_h(Y^0)$ for a fixed $\gamma>1$, we restart from $Y^0$ with an Armijo line search (backtracking factor $\beta_{\mathrm A}\in(0,1)$, sufficient-decrease parameter $c_{\mathrm A}\in(0,1)$, minimal step $a_{\min}>0$) for
\[
\Theta_h(Y)=\tfrac12\Delta t\,h^d\|\mathcal R_h(Y)\|_2^2.
\]
Since the interpolation weights are only piecewise smooth functions of the characteristic feet, Theorem~\ref{thm:reduced} does not establish quadratic convergence of the fully discrete method.

\section{Numerical experiments}\label{sec:numerics}\label{sec6:numerical}

We present four tests: a one-dimensional alignment problem, a smooth two-dimensional MFGC, an exact time-dependent benchmark, and a crowd-motion experiment. Tests 1, 2 and 4 examine the behavior of the Newton solvers, while Test 3 also measures discretization errors and compares reduced Newton--SL with reduced FD--Newton, where the fully nonlinear MFGC system is first discretized as in \cite{AK21}. We report residuals of the full discrete MFGC system, and check mass conservation and positivity. Discretization errors are assessed independently in Test 3 using the analytical solution. We give the Newton and GMRES parameters with the data of each test. GMRES settings are listed as (relative tolerance, restart size, maximum cycles). In the HJB march \eqref{eq:reducedsl-hjb}, the fixed-point iteration at each time level stops when the maximum-norm increment is below $10^{-13}$ ($10^{-11}$ in Test~3), with at most $100$ iterations followed by at most $100$ iterations damped by $\frac12$; a time level which does not reach the tolerance invalidates the residual evaluation. The Chebyshev iteration of Lemma~\ref{lem:Zdiscrete} stops at the relative residual $10^{-12}$ in the maximum norm, with at most $500$ iterations. In \eqref{eq:reducedsl-jv}, $\sigma_0=10^{-7}$; with inner tolerances between $10^{-13}$ and $10^{-11}$, the resulting error did not limit the final reduced residuals, which reach $10^{-9}$ to $10^{-13}$. The globalization uses $\gamma=10$, $\beta_{\mathrm A}=\frac12$, $c_{\mathrm A}=10^{-4}$ and $a_{\min}=2^{-12}$. In Tests~1, 2 and 4, the Newton stopping condition is based on the maximum norm of the accepted increment; in Test~3, it is based on the reduced residual $\mathcal R_h$. The full-system residual is evaluated independently at the computed solution, rather than used as an additional increment-based stopping condition. 

The full residual is measured by
\begin{equation}\label{eq:sec6:full-residual}
R_\infty=\max\bigl\{\|r_U\|_\infty,\|r_M\|_\infty,
\|r_Z\|_\infty\bigr\},
\end{equation}
where, for the alignment model and $k<N_t$,
\[
\begin{aligned}
(r_U)^k&=\frac{U^k-B^kU^{k+1}}{\Delta t}
  -L(x,\alpha^k,Z^k)-f(x,M^k),\\
(r_M)^k&=\frac{M^{k+1}-(B^k)^\top M^k}{\Delta t},\\
(r_Z)^k&=(K_hM^k+\varepsilon)Z^k
 -K_h[M^k(\lambda Z^k-D_hU^k)].
\end{aligned}
\]
For Test~3, the HJB residual $(r_U)^k$ additionally contains
the term $-V(t_k,x_i)$ at each grid node $x_i$.
The control residual is also evaluated at $k=N_t$; the terminal
and initial data are imposed directly. For the FD scheme in Test~3,
the residuals are evaluated using its corresponding discrete equations.
The residuals are reported as diagnostics, not as the Newton stopping
criteria in increment-based tests.

\subsection{Test 1: full Newton--SL in an alignment model}
\label{sec6:test1}\label{sec:test1}

We consider the one-dimensional alignment model \eqref{eq:align} on $\mathbb T$, with
\[
 T=1,\quad \nu=0.02,\quad f(m)=m,\quad
 G(x)=\tfrac12\bigl(1-\cos(2\pi(x-0.6))\bigr).
\]
On the spatial grid $x_i=ih$, the initial density is
\begin{equation}\label{eq:sec6:m01}
 m_{0,i}=\frac{0.5+\exp\!\left(-\bigl((x_i-0.3)/0.1\bigr)^2\right)}
 {h\sum_j\left[0.5+\exp\!\left(-\bigl((x_j-0.3)/0.1\bigr)^2\right)\right]}.
\end{equation}
The discrete Gaussian interaction weights are defined by
\begin{equation}\label{eq:sec6:K1}
\begin{aligned}
 d_{\mathbb T}(x,y)&=\min_{n\in\mathbb Z}|x-y+n|,\\
 \widetilde K_{ij}&=\exp\!\left(-\frac{d_{\mathbb T}(x_i,x_j)^2}{2\sigma_K^2}\right),\qquad
 w_{ij}=\frac{\widetilde K_{ij}}{\sum_\ell\widetilde K_{i\ell}},\qquad \sigma_K=0.1.
\end{aligned}
\end{equation}
We set $(K_h\phi)_i=\sum_jw_{ij}\phi_j$. We use $h=1/32$, $N_t=363$ and $\varepsilon=10^{-3}$. The direct initial guesses are $U^k=G$, $M^k=m_0$, $Z^k=0$ for full Newton--SL and $\Gamma^k=f(m_0)$, $Z^k=0$ for reduced Newton--SL. The full and reduced Newton--SL solvers use accepted-increment tolerances $10^{-8}$ and $10^{-9}$, respectively, with GMRES settings $(10^{-9},500,1)$ and $(10^{-8},100,3)$.

Table~\ref{tab:sec6:bgs} shows the spectrum and iteration counts for the first full Newton linear system. Relaxed BGS converges for $\lambda\leq0$, whereas both BGS variants diverge for $\lambda=0.5$ and $0.9$. GMRES converges in all tested cases. $s_1$ in Table~\ref{tab:sec6:bgs} is an estimate of the dominant eigenvalue of the BGS error operator, so that $|s_1|$ approximates $\mu$ of Proposition~\ref{prop:bgs} when $\lambda=0$; the sweeps stop when the maximum-norm increment is below $10^{-10}$, with at most $2000$ sweeps, and GMRES uses $(10^{-10},500,1)$. At $\lambda=0$, $s_1=-1.27$, i.e., $\mu\approx1.27>1$, so that plain BGS diverges by Proposition~\ref{prop:bgs}(a); for this problem, the method of \cite{CZ} requires relaxation or GMRES.

\begin{table}[htb!]
\centering\small
\caption{Test 1: dominant BGS eigenvalue and inner iteration counts for the first full Newton system. ``div.'' denotes divergence.}\label{tab:sec6:bgs}
\begin{tabular}{@{}rrrrr@{}}\toprule
$\lambda$&$s_1$& BGS& relaxed BGS& GMRES\\\midrule
$-0.9$&$-2.74$&div.&74&29\\
$-0.5$&$-2.23$&div.&52&25\\
$0$&$-1.27$&div.&28&17\\
$0.5$&$2.22$&div.&div.&72\\
$0.9$&$23.2$&div.&div.&192\\\bottomrule
\end{tabular}
\end{table}

\begin{figure}[htb!]
\centering
\includegraphics[width=0.92\linewidth]{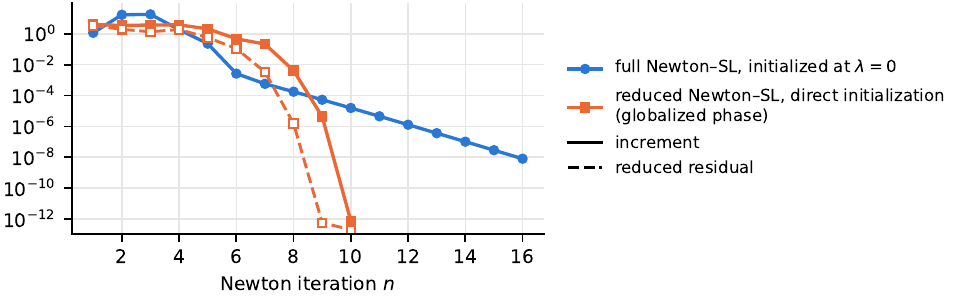}
\caption{Test 1, $\lambda=0.9$: nonlinear increment and reduced-residual histories.}
\label{fig:sec6:history1}
\end{figure}

At $\lambda=0.9$, full Newton--SL converges when initialized
with the solution for $\lambda=0$, but fails with the direct
guess after a Newton update produces a negative density.
Reduced Newton--SL converges from the direct guess and
retains a nonnegative density: along its iterations, the minimum of the reconstructed density is $1.2\cdot10^{-9}$ and the mass defect is below $10^{-15}$, whereas the failed full Newton--SL run reaches a density minimum of $-0.57$. Table~\ref{tab:sec6:test1out}
summarizes these results, and Figure~\ref{fig:sec6:history1} shows that the increments of full Newton--SL decrease linearly, whereas the reduced residual decreases quadratically in the final iterations. This is consistent with the construction of the two solvers: reduced Newton--SL differentiates the discrete residual $\mathcal R_h$, whereas full Newton--SL discretizes the continuous linearization \eqref{eq:newton}, so that its step, in which the derivative of $(B^k)^\top$ is replaced by $\mathcal N^k$, is not the exact Newton step of the discrete system. Figure~\ref{fig:sec6:history1} shows the globalized phase of the reduced run, which restarts from the direct guess after the two undamped steps counted in Table~\ref{tab:sec6:test1out}; there, ``local'' and ``global'' denote the undamped and globalized phases. The reduced run also needs this restart, and its cost is comparable to that of the failed full run.

\begin{table}[htb!]
\centering\small
\caption{Test 1, $\lambda=0.9$: Newton--SL outcomes. ``Full residual'' is the achieved value of $R_\infty$, independent of the increment-based termination criterion.}
\label{tab:sec6:test1out}
\begin{tabular}{@{}lllr r@{}}\toprule
method&initialization&Newton steps&CPU (s)&full residual\\\midrule
Full Newton--SL&from $\lambda=0$&16 local&68&$1.1\cdot10^{-7}$\\
Full Newton--SL&direct&failed: 2 local, 5 global&367&$1.2\cdot10^{2}$\\
Reduced Newton--SL&direct&2 local, 10 global&346&$2.6\cdot10^{-12}$\\\bottomrule
\end{tabular}
\end{table}

\begin{figure}[htb!]
\centering
\begin{minipage}[t]{0.49\linewidth}
\centering
\includegraphics[width=\linewidth]{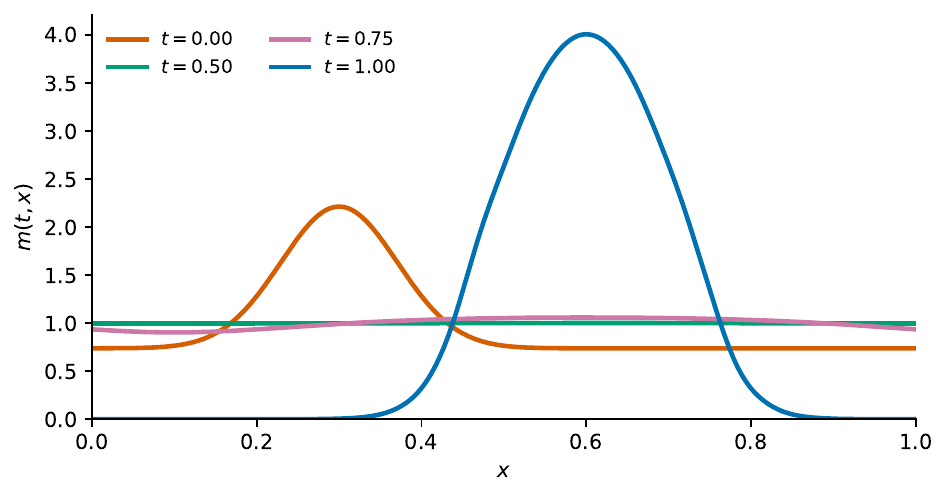}
\end{minipage}\hfill
\begin{minipage}[t]{0.49\linewidth}
\centering
\includegraphics[width=\linewidth]{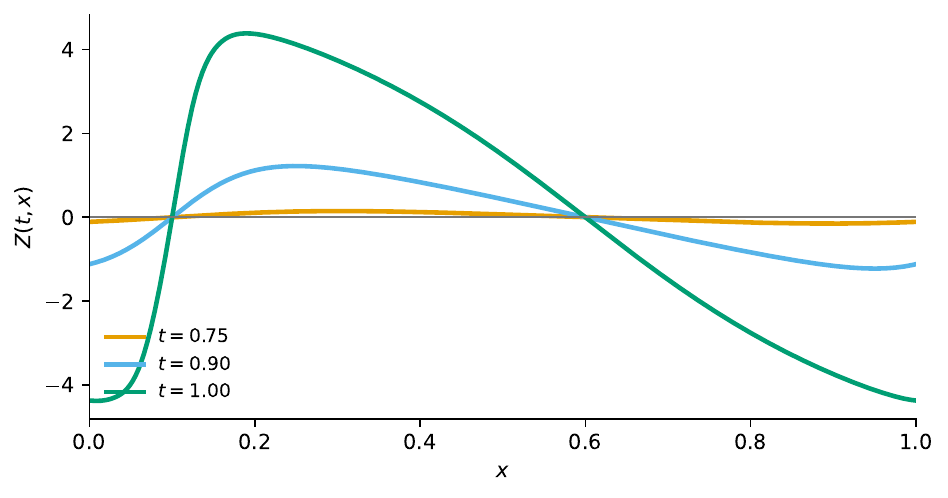}
\end{minipage}
\caption{Test 1, $\lambda=0.9$: density (left) and average-control (right) profiles for full Newton--SL.}
\label{fig:sec6:test1sol}
\end{figure}

Figure~\ref{fig:sec6:test1sol} displays the density and average control for the successful full Newton--SL run at $\lambda=0.9$. Initially concentrated near $x=0.3$, the density becomes almost uniform at $t\simeq T/2$ and then concentrates near $x=0.6$, where the terminal cost $G$ is minimized. The average control varies most strongly near the terminal time. 

\subsection{Test 2: a smooth, coupled two-dimensional problem}
\label{sec6:smooth}\label{sec:smooth2d}

We consider the alignment model on $\mathbb T^2$ with $T=0.5$, $\nu=0.02$, $\lambda=0.5$, $\varepsilon=10^{-3}$, and
\[
\begin{aligned}
 f(x,m)&=\tfrac12m,\\
 m_0(x)&=1+0.2\cos(2\pi x_1)\cos(2\pi x_2),\\
 G(x)&=0.15\bigl[2-\cos(2\pi(x_1-0.2))-\cos(2\pi(x_2+0.1))\bigr]
 +0.05\sin(2\pi(x_1+x_2)).
\end{aligned}
\]
The discrete interaction kernel is
\[
\begin{aligned}
 K(x_i,x_j)&=\frac{\widetilde K_{ij}}{h^2\sum_r\widetilde K_{ir}},\quad \widetilde K_{ij}=\exp\!\left(4\sum_{\ell=1}^2
 \bigl[\cos(2\pi(x_{i,\ell}-x_{j,\ell}))-1\bigr]\right).
\end{aligned}
\]
We use $N\times N$ spatial nodes with $N\in\{16,24,32,48,64\}$ and compare the time-step rules $\Delta t=h^{3/2}/2$ and $\Delta t= h/4$. Here $N_t=\lceil 2T/h^{3/2}\rceil$ and $N_t=\lceil 4T/h\rceil$, respectively, and $\Delta t=T/N_t$. Reduced Newton--SL starts from $\Gamma^k=f(m_0)$, $Z^k=0$, and uses the accepted-increment tolerance $10^{-9}$ and GMRES settings $(10^{-8},60,4)$. Table~\ref{tab:sec6:smooth} reports Newton iterations and CPU times. The method converges in five or six iterations in every case, with full-system residuals below $1.4\cdot10^{-11}$.

\begin{table}[htb!]
\centering\small
\caption{Test 2: reduced Newton--SL iteration counts $n_N$ and CPU times under two time-step rules.}
\label{tab:sec6:smooth}
\begin{tabular}{@{}r|rrr|rrr@{}}\toprule
&\multicolumn{3}{c|}{$\Delta t=h^{3/2}/2$}&
 \multicolumn{3}{c}{$\Delta t= h/4$}\\
$N$&$N_t$&$n_N$&CPU (s)&$N_t$&$n_N$&CPU (s)\\\midrule
16&64&5&33&32&5&25\\
24&118&5&94&48&5&58\\
32&182&5&195&64&6&131\\
48&333&5&701&96&6&675\\
64&512&5&1883&128&6&993\\\bottomrule
\end{tabular}
\end{table}

\begin{figure}[htb!]
\centering
\includegraphics[width=0.98\linewidth]{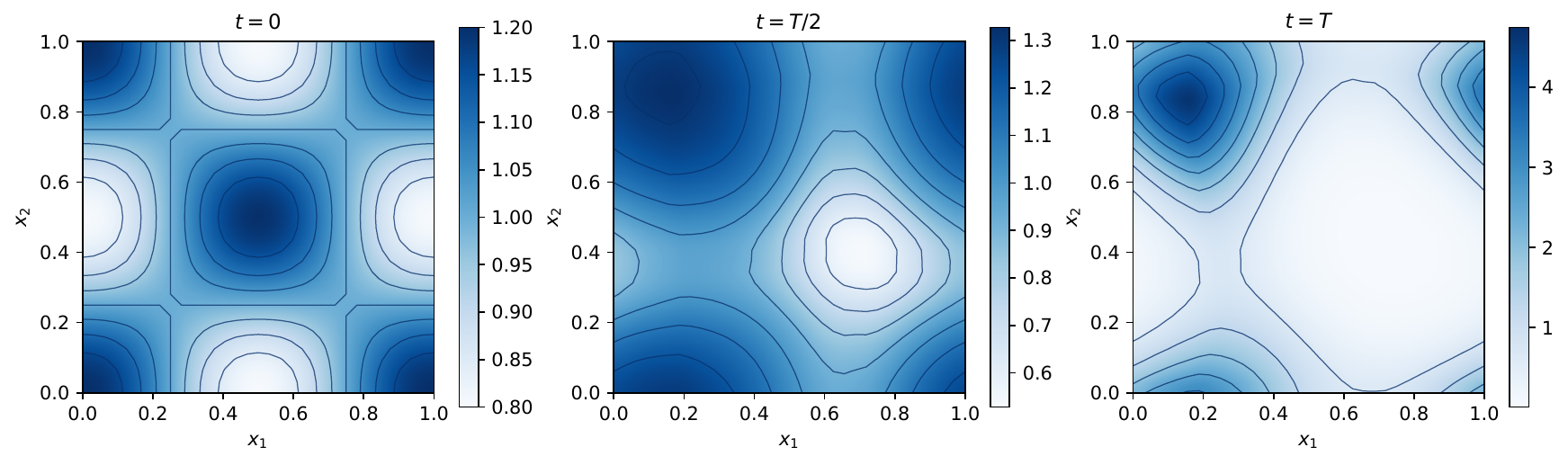}
\caption{Test 2: density at $t=0$, $T/2$ and $T$, with $N=32$ and $\Delta t=h^{3/2}/2$.}\label{fig:sec6:smooth-density}
\end{figure}

Figure~\ref{fig:sec6:smooth-density} shows the density evolution for $N=32$. The agents move toward the minimum of $G$, near which the density concentrates at $t=T$. Test 2 focuses on the nonlinear convergence of the solver; errors with respect to an analytical solution are studied in Test~3.

\subsection[Test 3: an exact, time-dependent MFGC benchmark and SL--FD accuracy--cost comparison]{Test 3: an exact, time-dependent MFGC benchmark and\\ SL--FD accuracy--cost comparison}
\label{sec6:exact}

The preceding experiments concern nonlinear convergence and do not measure discretization errors. We now construct an exact, genuinely time-dependent coupled MFGC solution on $[0,T]\times\mathbb T^2$, so that the errors in the value function, density and average control can be evaluated independently:

\begin{align}
 -\partial_t u-\nu\Delta u
  +\tfrac12|\nabla u|^2-\lambda Z\cdot\nabla u
      &=m^2+V(t,x),                                  \label{eq:sec6:hjb}\\
 \partial_t m-\nu\Delta m
      +\nabla\cdot\bigl[m(\lambda Z-\nabla u)\bigr]&=0, \label{eq:sec6:fp}\\
 Z(t,x)&=\frac{K*\bigl[m(\lambda Z-\nabla u)\bigr](t,x)}
                {(K*m)(t,x)+\varepsilon}.                 \label{eq:sec6:z}
\end{align}
Here $f(m)=m^2$ and $V=V(t,x)$ is a prescribed running potential. We take
\[
 T=0.1,\qquad \nu=0.03,\qquad \lambda=0.7,\qquad
 \varepsilon=10^{-3},\qquad \kappa=8.
\]
The periodic kernel and its Fourier multipliers are

\begin{equation}\label{eq:sec6:kernel}
 K(x)=\frac{\exp\{\kappa(\cos(2\pi x_1)+\cos(2\pi x_2))\}}
                {I_0(\kappa)^2},
 \qquad
 \widehat K(k)=\frac{I_{|k_1|}(\kappa)I_{|k_2|}(\kappa)}
                       {I_0(\kappa)^2},
\end{equation}
where $I_n$ denotes the modified Bessel function of the first kind.
Set $q=-0.38$ and $c=\nu+q=-0.35$. We consider the density

\begin{align}
 m^*(t,x)=1
 &+0.35e^{q(2\pi)^2t}\cos(2\pi x_1)
 +0.15e^{q(8\pi^2)t}\sin(2\pi(x_1+x_2))\notag\\[-1mm]
 &+0.10e^{q(2\pi)^2t}\cos(2\pi x_2).       \label{eq:sec6:exact-m}
\end{align}
For this choice $m^*\ge0.4$. With $W^*=K*m^*+\varepsilon$, define

\begin{equation}\label{eq:sec6:exact-uZ}
 Z^*=c\nabla\log W^*,\qquad
 u^*=c\bigl(\lambda\log W^*-\log m^*\bigr),\qquad
 \alpha^*=\lambda Z^*-\nabla u^*=c\nabla\log m^*.
\end{equation}
We set $m_0=m^*(0,\cdot)$ and $G=u^*(T,\cdot)$. The running
potential $V$ is chosen so that $(u^*,m^*,Z^*)$ satisfies the HJB
equation exactly. The FP and average-control equations are also
satisfied without source terms; see Appendix~\ref{app:exactMFGC}. Figure~\ref{fig:sec6:exact-control} shows $m^*$ and $Z^*$ at three times.
The time-dependent potential places this benchmark outside the
time-independent framework of Section~\ref{sec:theory}.

\begin{figure}[htb!]
\centering
\includegraphics[width=0.98\linewidth]{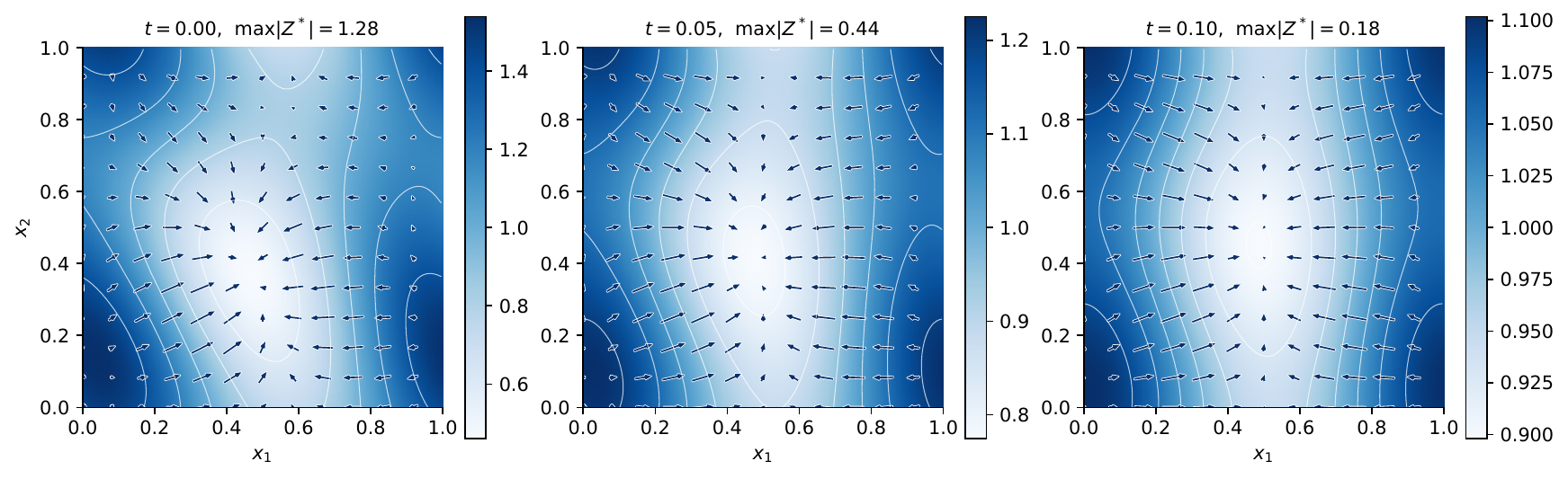}
\caption{Test 3: exact density and average control at three times.}
\label{fig:sec6:exact-control}
\end{figure}

We compare reduced Newton--SL with the reduced FD--Newton scheme based on the Engquist--Osher numerical Hamiltonian and adjoint FP discretization of \cite{AK21}. In the matched comparison, both methods use the same spatial grid, time step, initial and terminal data, and prescribed potential $V(t,x)$. The errors, evaluated against the analytical fields, are

\begin{equation}\label{eq:sec6:errors}
 E_u=\|U_h(0)-u^*(0)\|_{\ell_h^2},\qquad
 E_m=\|M_h(T)-m^*(T)\|_{\ell_h^2},\qquad
 E_Z=\|Z_h(0)-Z^*(0)\|_{\ell_h^2},
\end{equation}
where $\|e\|_{\ell_h^2}=(h^2\sum_i|e_i|^2)^{1/2}$. 

We compare both schemes with $N\in\{24,48,64,96\}$ and $N_t=N/2$. Both solvers start from $\Gamma^k=f(m_0)$, $Z^k=0$, and use the reduced-residual tolerance $5\cdot10^{-8}$ in the maximum norm and GMRES settings $(3\cdot10^{-7},36,2)$.
Table~\ref{tab:sec6:exactmatch} reports errors, CPU times, and the numbers of Newton iterations.

\begin{figure}[htb!]
\centering
\includegraphics[width=0.98\linewidth]{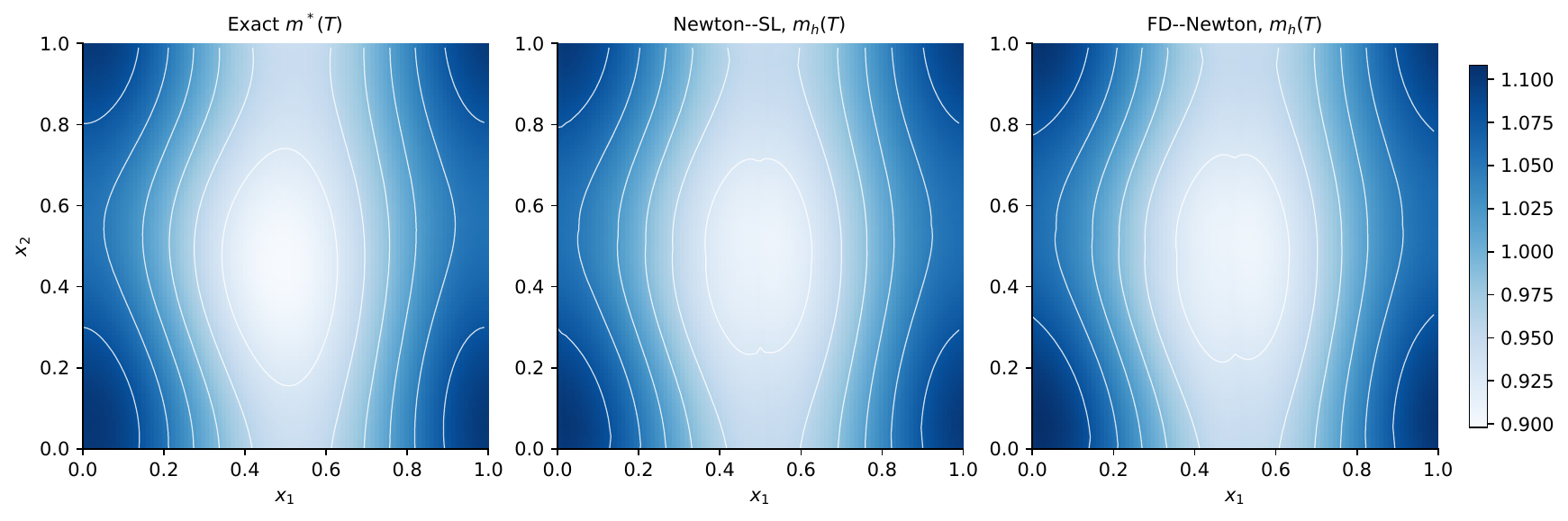}
\caption{Test 3, $N=96$, $N_t=48$: exact and computed densities at $t=T$.}\label{fig:sec6:exact-density}
\end{figure}

\begin{figure}[htb!]
\centering
\includegraphics[width=0.68\linewidth]{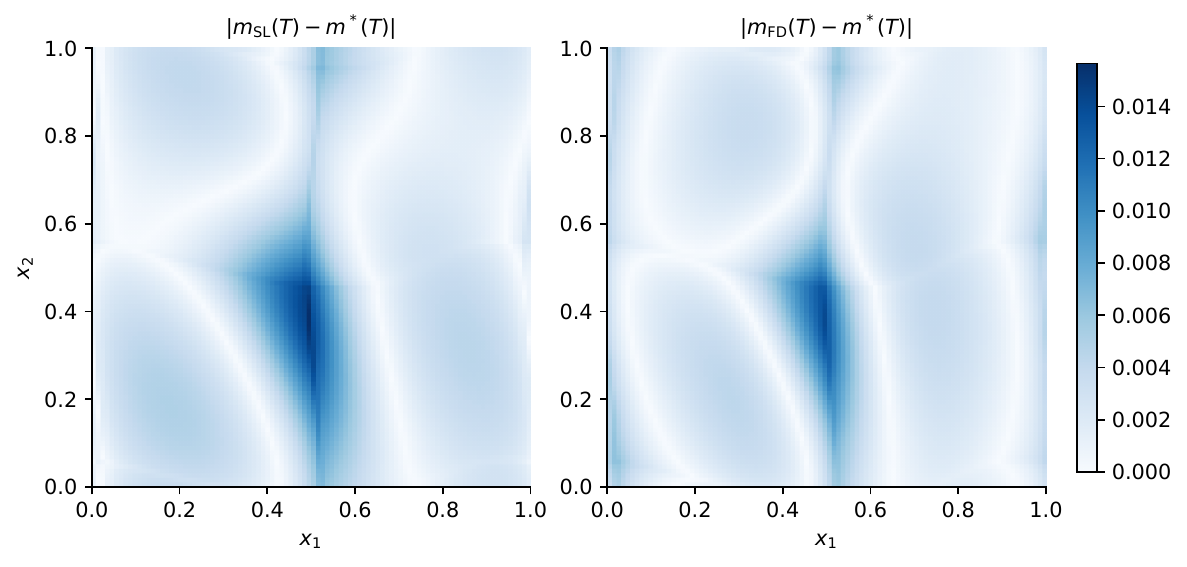}
\caption{Test 3, $N=96$, $N_t=48$: pointwise absolute errors of the terminal density for SL and FD.}\label{fig:sec6:exact-error}
\end{figure}

Both schemes converge at about first order, which is the rate expected for $\Delta t\propto h$ from the formal SL error $O(\Delta t+h^2/\Delta t)$. When $h$ is divided by four, from $N=24$ to $N=96$, the SL errors $E_u$, $E_m$ and $E_Z$ decrease by the factors $5.0$, $3.7$ and $4.0$, and the FD errors by the factors $3.7$, $3.8$ and $4.3$. SL has the smaller errors in $u$ and $Z$, by a factor between $1.7$ and $2.5$, and the errors in $m$ are comparable, with FD slightly smaller. The numbers of Newton iterations are similar. SL is faster on every grid, since one FD residual evaluation requires a Newton iteration with a sparse factorization at each time level. Figures~\ref{fig:sec6:exact-density} and~\ref{fig:sec6:exact-error} show the densities at $t=T$ and the pointwise errors for $N=96$.

\begin{table}[htb!]
\centering\small
\caption{Test 3: errors against the exact solution, CPU times, and numbers of Newton iterations for SL and FD, with $N_t=N/2$. CPU times on one core of an Intel Xeon processor at 2.1 GHz.}
\label{tab:sec6:exactmatch}
\begin{tabular}{@{}rrlrrrrr@{}}\toprule
$N$&$N_t$&scheme&CPU (s)&Newton&$E_u$&$E_m$&$E_Z$\\\midrule
24&12&SL&3&3&$1.18\cdot10^{-3}$&$1.44\cdot10^{-2}$&$1.57\cdot10^{-2}$\\
24&12&FD&16&5&$2.16\cdot10^{-3}$&$1.20\cdot10^{-2}$&$3.98\cdot10^{-2}$\\
48&24&SL&20&4&$6.54\cdot10^{-4}$&$7.05\cdot10^{-3}$&$9.86\cdot10^{-3}$\\
48&24&FD&102&5&$1.15\cdot10^{-3}$&$6.23\cdot10^{-3}$&$1.92\cdot10^{-2}$\\
64&32&SL&33&3&$3.95\cdot10^{-4}$&$5.45\cdot10^{-3}$&$5.99\cdot10^{-3}$\\
64&32&FD&193&4&$8.72\cdot10^{-4}$&$4.72\cdot10^{-3}$&$1.42\cdot10^{-2}$\\
96&48&SL&127&4&$2.38\cdot10^{-4}$&$3.87\cdot10^{-3}$&$3.91\cdot10^{-3}$\\
96&48&FD&971&4&$5.90\cdot10^{-4}$&$3.18\cdot10^{-3}$&$9.37\cdot10^{-3}$\\\bottomrule
\end{tabular}
\end{table}

In all runs, the componentwise full-system residuals are below $8.0\cdot10^{-8}$, considerably below the discretization errors, the computed density stays above $0.46$, and the mass defect is below $10^{-14}$. Stability of the computed solutions in the sense of Definition~\ref{def:stable} is not verified; the quadratic decrease of the reduced residual in the final iterations, e.g.\ $6.2\cdot10^{-3}$, $1.1\cdot10^{-5}$, $5.8\cdot10^{-13}$ for SL with $N=48$, is consistent with a nonsingular Jacobian of $\mathcal R_h$, the discrete counterpart of the characterization in Theorem~\ref{thm:reduced}(b).

\begin{figure}[htb!]
\centering
\includegraphics[width=0.98\linewidth]{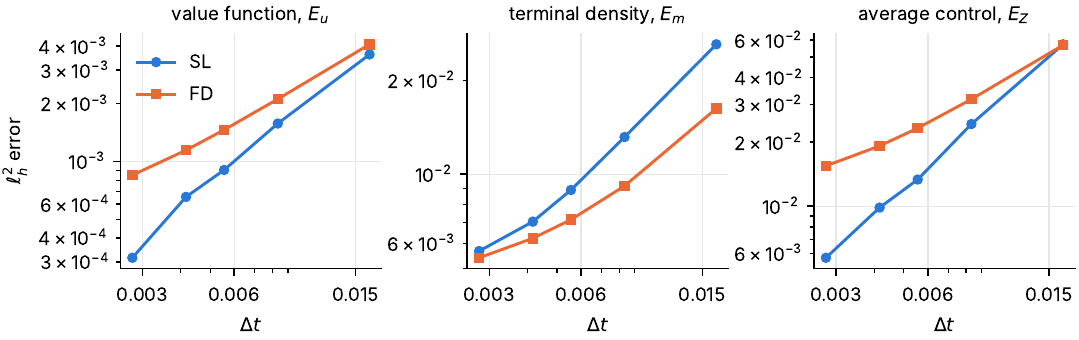}
\caption{Test 3: time-step dependence of the three errors at $N=48$.}
\label{fig:sec6:timesteps}
\end{figure}

We next vary $N_t\in\{6,12,18,24,36\}$ at fixed $N=48$. Figure~\ref{fig:sec6:timesteps} shows the errors as functions of $\Delta t$. The errors of both schemes decrease monotonically with $\Delta t$.

\FloatBarrier
\subsection{Test 4: strong alignment and small viscosity in crowd motion}
\label{sec6:crowd}\label{sec:ak2d}

We consider a two-dimensional crowd-motion problem inspired by \cite[Section~4.1]{AK21}. On $\Omega=(-1/2,1/2)^2$, with reflecting boundary conditions, the data are
\[
 T=1,\qquad \nu=10^{-3},\qquad \lambda=0.9,\qquad
 f(x,m)=10^{-3}m+0.1\varphi(x),\qquad G(x)=\varphi(x).
\]
With the smoothed indicator
\[
\chi_{a}(x;c)=\tfrac14\Big(1+\tanh\frac{a-|x_1-c_1|}{\delta}\Big)\Big(1+\tanh\frac{a-|x_2-c_2|}{\delta}\Big),\qquad \delta=0.02,
\]
of the square of half-width $a$ centered at $c\in\Omega$, the initial density and the potential are
\begin{align*}
m_0(x)&=C_0\Big(\chi_{0.1}\big(x;(0.3,0.3)\big)+\chi_{0.1}\big(x;(-0.3,-0.3)\big)+10^{-3}\Big),\\
\varphi(x)&=1-\chi_{0.15}\big(x;(-0.35,0.35)\big)-\chi_{0.15}\big(x;(0.35,-0.35)\big),
\end{align*}
where $C_0>0$ is such that $\int_\Omega m_0=1$. The interaction kernel is $K(x,y)=\kappa(|x-y|)$ with $\rho=0.2$ and
\[
\kappa(r)=\begin{cases}1,& r\le0.9\rho,\\ 1-3s^2+2s^3,\quad s=\dfrac{r-0.9\rho}{0.1\rho},& 0.9\rho<r<\rho,\\ 0,& r\ge\rho.\end{cases}
\]
These data differ from those in \cite{AK21}.
 We take $\varepsilon=10^{-3}$ and $N_t=\lceil 2/h^{3/2}\rceil$. Reduced Newton--SL uses the accepted-increment tolerance $10^{-7}$ and GMRES settings $(10^{-6},200,1)$; the direct initialization is $\Gamma^k=f(m_0)$, $Z^k=0$.

\begin{figure}[htb!]
\centering
\includegraphics[width=\linewidth]{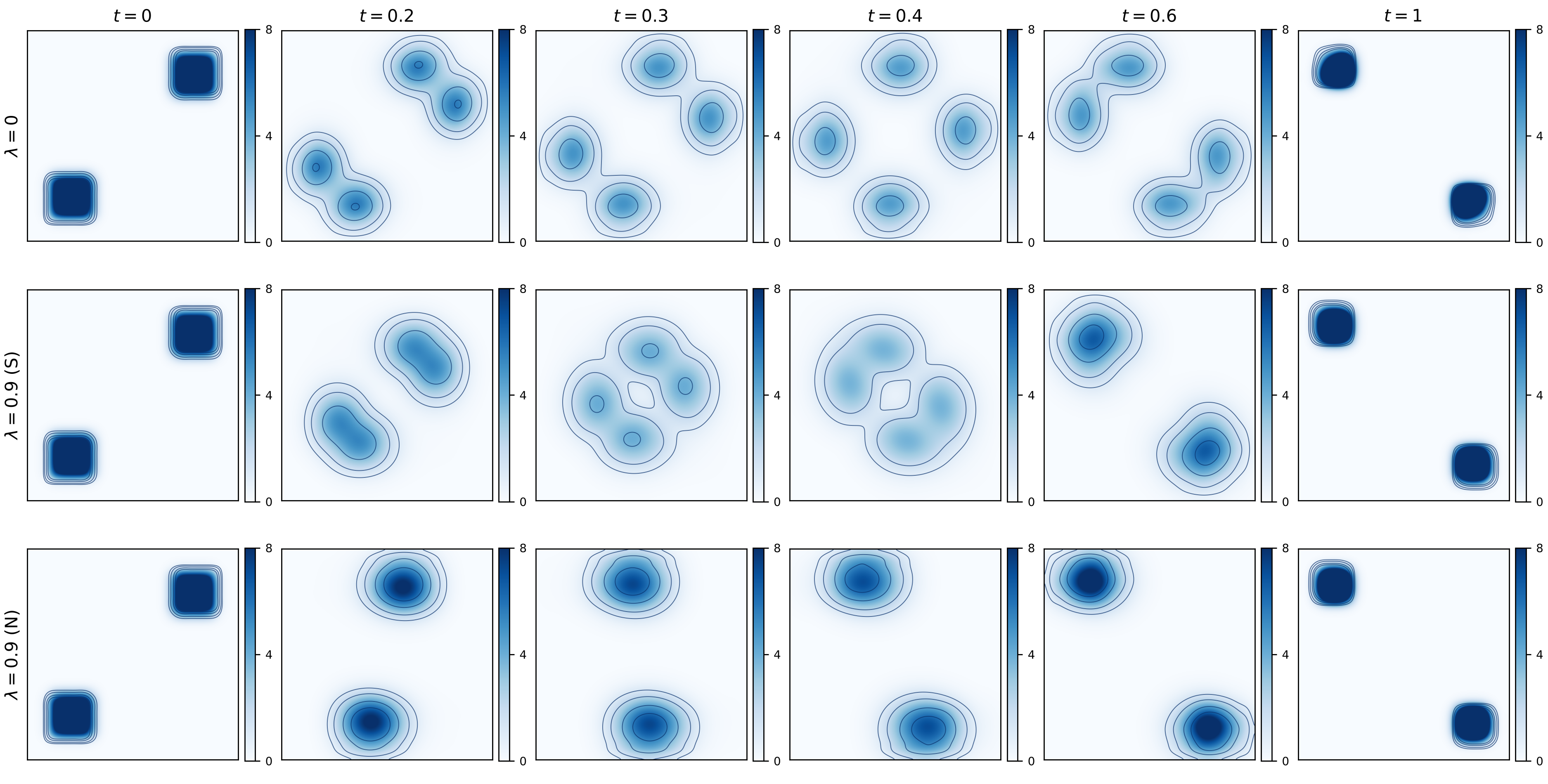}
\caption{Test 4, $h=1/26$: density snapshots at $\lambda=0$ and at $\lambda=0.9$ for the symmetric (S) and nonsymmetric (N) discrete branches.}
\label{fig:sec6:crowd}
\end{figure}

Figure~\ref{fig:sec6:crowd} shows two distinct configurations at $\lambda=0.9$. On the symmetric branch (S), the two crowds separate and recombine before approaching their targets, whereas on the nonsymmetric branch (N) they move predominantly toward the targets along more direct paths. The symmetric runs are initialized with the solution for $\lambda=0$; the nonsymmetric configuration is obtained directly on the coarser grid and then used to initialize the finer grid. The iteration counts and the achieved full residuals are reported in Table~\ref{tab:sec6:crowd}; in all four runs the minimum of the reconstructed density is nonnegative to rounding and the mass defect is below $10^{-15}$. These configurations solve the discrete SL system and are qualitatively consistent with the patterns in \cite{AK21}.

\begin{table}[htb!]
\centering\small
\caption{Test 4: reduced Newton--SL iterations and achieved full residuals for the two discrete branches at $\lambda=0.9$.}
\label{tab:sec6:crowd}
\begin{tabular}{@{}lllr l r@{}}\toprule
$h$ & $N_t$ & initialization & Newton iterations & branch & $R_\infty$\\\midrule
$1/26$ & 266 & $\lambda=0$ solution & 6 & S & $2.5\cdot10^{-11}$\\
$1/26$ & 266 & direct & 10 & N & $2.1\cdot10^{-11}$\\
$1/51$ & 729 & $\lambda=0$ solution & 6 & S & $7.6\cdot10^{-11}$\\
$1/51$ & 729 & refined N branch & 4 & N & $8.5\cdot10^{-11}$\\\bottomrule
\end{tabular}
\end{table}

\FloatBarrier
\section{Conclusions}\label{sec:conclusions}

We proved the local quadratic convergence of full and reduced Newton iterations for second-order MFGC near stable solutions, and proposed an SL scheme with two matrix-free Newton solvers. In the experiments, the reduced Newton--GMRES solver, which differentiates the discrete residual and keeps the density nonnegative and of unit mass, converges quadratically in its final iterations, whereas full Newton--SL, which discretizes the continuous linearization, converges linearly and can fail without a suitable initialization. On the exact benchmark, the SL and finite-difference schemes both converge at first order with errors of comparable size, and the SL implementation is faster.

Several questions remain open: the convergence of the fully discrete SL scheme and of the discrete Newton iterations, the solvability of the implicit HJB step, the stability of the computed solutions, in particular of the crowd-motion branches, and the extension of the analysis to bounded domains with reflecting boundary conditions. The exact solution of Appendix~\ref{app:exactMFGC} can be used as a benchmark for other MFGC discretizations.

\appendix
\section{Construction of the exact time-dependent coupled MFGC benchmark}\label{app:exactMFGC}
We derive the analytical solution used in Test~3 and the corresponding running potential. The FP and average-control equations are satisfied without source terms.

\paragraph{A positive density solving a scalar linear equation}
For a constant $q<0$, choose finitely many nonzero integer wave vectors $k_j\in\mathbb Z^2$ and mean-zero trigonometric modes $\phi_j$ satisfying $\Delta\phi_j=-4\pi^2|k_j|^2\phi_j$. Set
\begin{equation}\label{eq:app:m}
 m^*(t,x)=1+\sum_{j=1}^J a_j e^{4\pi^2 q|k_j|^2t}\phi_j(x),\qquad
 \sum_{j=1}^J|a_j|<1,
\end{equation}
where $\|\phi_j\|_\infty\le1$. Hence $m^*>0$, $\int_{\mathbb T^2}m^*(t,x)\,\mathrm dx=1$, and differentiation gives
\begin{equation}\label{eq:app:heat}
 \partial_t m^*+q\Delta m^*=0.
\end{equation}
The three modes and coefficients in~\eqref{eq:sec6:exact-m} correspond to $k_1=(1,0)$, $k_2=(1,1)$, $k_3=(0,1)$, with amplitudes $0.35$, $0.15$, and $0.10$. Since $q<0$, all amplitudes decrease in magnitude for $t\ge0$, and $m^*\ge0.4$ for these choices.

\paragraph{Nonlocal control identity}
Let $K$ denote the smooth periodic kernel in~\eqref{eq:sec6:kernel}.
It is nonnegative and translation invariant. Set
\[
 W^*=K*m^*+\varepsilon,
 \qquad c=\nu+q,\qquad
 \alpha^*=c\nabla\log m^*,\qquad
 Z^*=c\nabla\log W^*.
\]
Since $K$ is independent of $t$ and is a periodic convolution kernel, spatial differentiation commutes with convolution. Consequently,
\[
 K*(m^*\alpha^*)=cK*\nabla m^*
   =c\nabla(K*m^*)=c\nabla W^*=W^* Z^*.
\]
Thus the averaged-control relation~\eqref{eq:sec6:z} holds pointwise. Define
\[
 u^*=c\bigl(\lambda\log W^*-\log m^*\bigr).
\]
Then $\lambda Z^*-\nabla u^*=c\nabla\log m^*=\alpha^*$, and
\begin{align*}
 \partial_t m^*-\nu\Delta m^*+\operatorname{div}(m^*\alpha^*)
 &=\partial_t m^*-\nu\Delta m^*+c\Delta m^*\\
 &=\partial_t m^*+q\Delta m^*=0
\end{align*}
by~\eqref{eq:app:heat}. This identity is exact, with no source added to the FP equation.

\paragraph{Exact running potential}
Write $r=\log m^*$ and $s=\log W^*$. For the alignment Hamiltonian $H(p,Z)=|p|^2/2-\lambda Z\cdot p$ and the coupling $m^2$, the prescribed potential is
\begin{equation}\label{eq:app:Vcompact}
 V=-c(\lambda s_t-r_t)-\nu c(\lambda\Delta s-\Delta r)
   +\frac{c^2}{2}\bigl(|\nabla r|^2-\lambda^2|\nabla s|^2\bigr)
   -(m^*)^2.
\end{equation}
Indeed $\nabla u^*=c(\lambda\nabla s-\nabla r)$ and $Z^*=c\nabla s$; the cross terms cancel in
\[
 \tfrac12|\nabla u^*|^2-\lambda Z^*\cdot\nabla u^*
 =\tfrac{c^2}{2}\bigl(|\nabla r|^2-\lambda^2|\nabla s|^2\bigr).
\]
Moreover $\partial_tW^*=K*\partial_tm^*=-q\Delta W^*$, so $r_t=-q\Delta m^*/m^*$ and $s_t=-q\Delta W^*/W^*$. Together with the elementary logarithmic identities
\[
 \Delta\log a=\frac{\Delta a}{a}-\frac{|\nabla a|^2}{a^2},
 \qquad \nabla\log a=\frac{\nabla a}{a},
\]
formula~\eqref{eq:app:Vcompact} can be evaluated using only $m^*$, $W^*$ and their spatial derivatives. Substitution into the HJB equation~\eqref{eq:sec6:hjb} yields an identity by construction.

\paragraph{Analytical convolution and data}
The von Mises kernel is normalized by the identity
\[
\int_0^1 e^{\kappa\cos(2\pi y)}\,\mathrm dy=I_0(\kappa).
\]
For each Fourier wave vector $k=(k_1,k_2)$, its multiplier is
\[
 \widehat K(k)=\frac{I_{|k_1|}(\kappa)I_{|k_2|}(\kappa)}{I_0(\kappa)^2}.
\]
Therefore the exact $K*m^*$ is obtained by multiplying each amplitude in~\eqref{eq:sec6:exact-m} by the corresponding $\widehat K(k_j)$; no grid-based convolution is needed. Finally we impose $m(0)=m^*(0)$ and $u(T)=u^*(T)$, so there is no incompatibility at the initial or terminal time. For the parameters of Test~3, $c=-0.35\ne0$ and the spatially nonconstant $W^*$ gives $Z^*\not\equiv0$, ensuring that the benchmark genuinely exercises interaction through controls.
\section*{Declarations}






\subsection*{Use of AI tools}
The author used Claude (Anthropic) to assist with language editing, LaTeX formatting, code assisting and  manuscript proofreading. All results, proofs, and numerical experiments were developed and verified by the author, who takes full responsibility for the content of this article.

\end{document}